\newif\ifams
\newif\ifmscs
\amstrue
\ifams
\documentclass[reqno,11pt]{amsart}
\usepackage[foot]{amsaddr}
\fi
\ifmscs
\documentclass[11pt]{mscs}
\fi

\usepackage{tikz}
\usepackage{tikz-cd}
\usetikzlibrary{calc}
\usepackage{amsmath}
\usepackage{amssymb}
\usepackage{amstext}
\usepackage{amsfonts}
\usepackage{xspace}
\usepackage{url}
\usepackage{enumerate}
\long\def\ig#1{\relax}
\ig{Thanks to Roberto Minio for this def'n.  Compare the def'n of
\comment in AMSTeX.}
 
\newcount \coefa
\newcount \coefb
\newcount \coefc
\newcount\tempcounta
\newcount\tempcountb
\newcount\tempcountc
\newcount\tempcountd
\newcount\xext
\newcount\yext
\newcount\xoff
\newcount\yoff
\newcount\gap%
\newcount\arrowtypea
\newcount\arrowtypeb
\newcount\arrowtypec
\newcount\arrowtyped
\newcount\arrowtypee
\newcount\height
\newcount\width
\newcount\xpos
\newcount\ypos
\newcount\run
\newcount\rise
\newcount\arrowlength
\newcount\halflength
\newcount\arrowtype
\newdimen\tempdimen
\newdimen\xlen
\newdimen\ylen
\newsavebox{\tempboxa}%
\newsavebox{\tempboxb}%
\newsavebox{\tempboxc}%
 
\makeatletter
\def\settypes(#1,#2,#3){\arrowtypea#1 \arrowtypeb#2 \arrowtypec#3}
\def\settoheight#1#2{\setbox\@tempboxa\hbox{#2}#1\ht\@tempboxa\relax}%
\def\settodepth#1#2{\setbox\@tempboxa\hbox{#2}#1\dp\@tempboxa\relax}%
\def\settokens[#1`#2`#3`#4]{%
     \def\tokena{#1}\def\tokenb{#2}\def\tokenc{#3}\def\tokend{#4}}
\def\setsqparms[#1`#2`#3`#4;#5`#6]{%
\arrowtypea #1
\arrowtypeb #2
\arrowtypec #3
\arrowtyped #4
\width #5
\height #6
}
\def\setpos(#1,#2){\xpos=#1 \ypos#2}
 
\def\bfig{\begin{picture}(\xext,\yext)(\xoff,\yoff)}
\def\efig{\end{picture}}
 
\def\putbox(#1,#2)#3{\put(#1,#2){\makebox(0,0){$#3$}}}
 
\def\settriparms[#1`#2`#3;#4]{\settripairparms[#1`#2`#3`1`1;#4]}%
 
\def\settripairparms[#1`#2`#3`#4`#5;#6]{%
\arrowtypea #1
\arrowtypeb #2
\arrowtypec #3
\arrowtyped #4
\arrowtypee #5
\width #6
\height #6
}
 
\def\resetparms{\settripairparms[1`1`1`1`1;500]\width 500}
 
\resetparms
 
\def\mvector(#1,#2)#3{
\put(0,0){\vector(#1,#2){#3}}%
\put(0,0){\vector(#1,#2){30}}%
}
\def\evector(#1,#2)#3{{
\arrowlength #3
\put(0,0){\vector(#1,#2){\arrowlength}}%
\advance \arrowlength by-30
\put(0,0){\vector(#1,#2){\arrowlength}}%
}}
 
\def\horsize#1#2{%
\settowidth{\tempdimen}{$#2$}%
#1=\tempdimen
\divide #1 by\unitlength
}
 
\def\vertsize#1#2{%
\settoheight{\tempdimen}{$#2$}%
#1=\tempdimen
\settodepth{\tempdimen}{$#2$}%
\advance #1 by\tempdimen
\divide #1 by\unitlength
}
 
\def\vertadjust[#1`#2`#3]{%
\vertsize{\tempcounta}{#1}%
\vertsize{\tempcountb}{#2}%
\ifnum \tempcounta<\tempcountb \tempcounta=\tempcountb \fi
\divide\tempcounta by2
\vertsize{\tempcountb}{#3}%
\ifnum \tempcountb>0 \advance \tempcountb by20 \fi
\ifnum \tempcounta<\tempcountb \tempcounta=\tempcountb \fi
}
 
\def\horadjust[#1`#2`#3]{%
\horsize{\tempcounta}{#1}%
\horsize{\tempcountb}{#2}%
\ifnum \tempcounta<\tempcountb \tempcounta=\tempcountb \fi
\divide\tempcounta by20
\horsize{\tempcountb}{#3}%
\ifnum \tempcountb>0 \advance \tempcountb by60 \fi
\ifnum \tempcounta<\tempcountb \tempcounta=\tempcountb \fi
}
 
\ig{ In this procedure, #1 is the paramater that sticks out all the way,
#2 sticks out the least and #3 is a label sticking out half way.  #4 is
the amount of the offset.}
 
\def\sladjust[#1`#2`#3]#4{%
\tempcountc=#4
\horsize{\tempcounta}{#1}%
\divide \tempcounta by2
\horsize{\tempcountb}{#2}%
\divide \tempcountb by2
\advance \tempcountb by-\tempcountc
\ifnum \tempcounta<\tempcountb \tempcounta=\tempcountb\fi
\divide \tempcountc by2
\horsize{\tempcountb}{#3}%
\advance \tempcountb by-\tempcountc
\ifnum \tempcountb>0 \advance \tempcountb by80\fi
\ifnum \tempcounta<\tempcountb \tempcounta=\tempcountb\fi
\advance\tempcounta by20
}
 
\def\putvector(#1,#2)(#3,#4)#5#6{{%
\xpos=#1
\ypos=#2
\run=#3
\rise=#4
\arrowlength=#5
\arrowtype=#6
\ifnum \arrowtype<0
    \ifnum \run=0
        \advance \ypos by-\arrowlength
    \else
        \tempcounta \arrowlength
        \multiply \tempcounta by\rise
        \divide \tempcounta by\run
        \ifnum\run>0
            \advance \xpos by\arrowlength
            \advance \ypos by\tempcounta
        \else
            \advance \xpos by-\arrowlength
            \advance \ypos by-\tempcounta
        \fi
    \fi
    \multiply \arrowtype by-1
    \multiply \rise by-1
    \multiply \run by-1
\fi
\ifnum \arrowtype=1
    \put(\xpos,\ypos){\vector(\run,\rise){\arrowlength}}%
\else\ifnum \arrowtype=2
    \put(\xpos,\ypos){\mvector(\run,\rise)\arrowlength}%
\else\ifnum\arrowtype=3
    \put(\xpos,\ypos){\evector(\run,\rise){\arrowlength}}%
\fi\fi\fi
}}
 
\def\putsplitvector(#1,#2)#3#4{
\xpos #1
\ypos #2
\arrowtype #4
\halflength #3
\arrowlength #3
\gap 140
\advance \halflength by-\gap
\divide \halflength by2
\ifnum \arrowtype=1
    \put(\xpos,\ypos){\line(0,-1){\halflength}}%
    \advance\ypos by-\halflength
    \advance\ypos by-\gap
    \put(\xpos,\ypos){\vector(0,-1){\halflength}}%
\else\ifnum \arrowtype=2
    \put(\xpos,\ypos){\line(0,-1)\halflength}%
    \put(\xpos,\ypos){\vector(0,-1)3}%
    \advance\ypos by-\halflength
    \advance\ypos by-\gap
    \put(\xpos,\ypos){\vector(0,-1){\halflength}}%
\else\ifnum\arrowtype=3
    \put(\xpos,\ypos){\line(0,-1)\halflength}%
    \advance\ypos by-\halflength
    \advance\ypos by-\gap
    \put(\xpos,\ypos){\evector(0,-1){\halflength}}%
\else\ifnum \arrowtype=-1
    \advance \ypos by-\arrowlength
    \put(\xpos,\ypos){\line(0,1){\halflength}}%
    \advance\ypos by\halflength
    \advance\ypos by\gap
    \put(\xpos,\ypos){\vector(0,1){\halflength}}%
\else\ifnum \arrowtype=-2
    \advance \ypos by-\arrowlength
    \put(\xpos,\ypos){\line(0,1)\halflength}%
    \put(\xpos,\ypos){\vector(0,1)3}%
    \advance\ypos by\halflength
    \advance\ypos by\gap
    \put(\xpos,\ypos){\vector(0,1){\halflength}}%
\else\ifnum\arrowtype=-3
    \advance \ypos by-\arrowlength
    \put(\xpos,\ypos){\line(0,1)\halflength}%
    \advance\ypos by\halflength
    \advance\ypos by\gap
    \put(\xpos,\ypos){\evector(0,1){\halflength}}%
\fi\fi\fi\fi\fi\fi
}
 
\def\putmorphism(#1)(#2,#3)[#4`#5`#6]#7#8#9{{%
\run #2
\rise #3
\ifnum\rise=0
  \puthmorphism(#1)[#4`#5`#6]{#7}{#8}{#9}%
\else\ifnum\run=0
  \putvmorphism(#1)[#4`#5`#6]{#7}{#8}{#9}%
\else
\setpos(#1)%
\arrowlength #7
\arrowtype #8
\ifnum\run=0
\else\ifnum\rise=0
\else
\ifnum\run>0
    \coefa=1
\else
   \coefa=-1
\fi
\ifnum\arrowtype>0
   \coefb=0
   \coefc=-1
\else
   \coefb=\coefa
   \coefc=1
   \arrowtype=-\arrowtype
\fi
\width=2
\multiply \width by\run
\divide \width by\rise
\ifnum \width<0  \width=-\width\fi
\advance\width by60
\if l#9 \width=-\width\fi
\putbox(\xpos,\ypos){#4}
{\multiply \coefa by\arrowlength
\advance\xpos by\coefa
\multiply \coefa by\rise
\divide \coefa by\run
\advance \ypos by\coefa
\putbox(\xpos,\ypos){#5} }%
{\multiply \coefa by\arrowlength
\divide \coefa by2
\advance \xpos by\coefa
\advance \xpos by\width
\multiply \coefa by\rise
\divide \coefa by\run
\advance \ypos by\coefa
\if l#9%
   \put(\xpos,\ypos){\makebox(0,0)[r]{$#6$}}%
\else\if r#9%
   \put(\xpos,\ypos){\makebox(0,0)[l]{$#6$}}%
\fi\fi }%
{\multiply \rise by-\coefc
\multiply \run by-\coefc
\multiply \coefb by\arrowlength
\advance \xpos by\coefb
\multiply \coefb by\rise
\divide \coefb by\run
\advance \ypos by\coefb
\multiply \coefc by70
\advance \ypos by\coefc
\multiply \coefc by\run
\divide \coefc by\rise
\advance \xpos by\coefc
\multiply \coefa by140
\multiply \coefa by\run
\divide \coefa by\rise
\advance \arrowlength by\coefa
\ifnum \arrowtype=1
   \put(\xpos,\ypos){\vector(\run,\rise){\arrowlength}}%
\else\ifnum\arrowtype=2
   \put(\xpos,\ypos){\mvector(\run,\rise){\arrowlength}}%
\else\ifnum\arrowtype=3
   \put(\xpos,\ypos){\evector(\run,\rise){\arrowlength}}%
\fi\fi\fi}\fi\fi\fi\fi}}
 
\def\puthmorphism(#1,#2)[#3`#4`#5]#6#7#8{{%
\xpos #1
\ypos #2
\width #6
\arrowlength #6
\putbox(\xpos,\ypos){#3\vphantom{#4}}%
{\advance \xpos by\arrowlength
\putbox(\xpos,\ypos){\vphantom{#3}#4}}%
\horsize{\tempcounta}{#3}%
\horsize{\tempcountb}{#4}%
\divide \tempcounta by2
\divide \tempcountb by2
\advance \tempcounta by30
\advance \tempcountb by30
\advance \xpos by\tempcounta
\advance \arrowlength by-\tempcounta
\advance \arrowlength by-\tempcountb
\putvector(\xpos,\ypos)(1,0){\arrowlength}{#7}%
\divide \arrowlength by2
\advance \xpos by\arrowlength
\vertsize{\tempcounta}{#5}%
\divide\tempcounta by2
\advance \tempcounta by20
\if a#8 %
   \advance \ypos by\tempcounta
   \putbox(\xpos,\ypos){#5}%
\else
   \advance \ypos by-\tempcounta
   \putbox(\xpos,\ypos){#5}%
\fi}}
 
\def\putvmorphism(#1,#2)[#3`#4`#5]#6#7#8{{%
\xpos #1
\ypos #2
\arrowlength #6
\arrowtype #7
\settowidth{\xlen}{$#5$}%
\putbox(\xpos,\ypos){#3}%
{\advance \ypos by-\arrowlength
\putbox(\xpos,\ypos){#4}}%
{\advance\arrowlength by-140
\advance \ypos by-70
\ifdim\xlen>0pt
   \if m#8%
      \putsplitvector(\xpos,\ypos){\arrowlength}{\arrowtype}%
   \else
      \putvector(\xpos,\ypos)(0,-1){\arrowlength}{\arrowtype}%
   \fi
\else
   \putvector(\xpos,\ypos)(0,-1){\arrowlength}{\arrowtype}%
\fi}%
\ifdim\xlen>0pt
   \divide \arrowlength by2
   \advance\ypos by-\arrowlength
   \if l#8%
      \advance \xpos by-40
      \put(\xpos,\ypos){\makebox(0,0)[r]{$#5$}}%
   \else\if r#8%
      \advance \xpos by40
      \put(\xpos,\ypos){\makebox(0,0)[l]{$#5$}}%
   \else
      \putbox(\xpos,\ypos){#5}%
   \fi\fi
\fi
}}
 
\def\topadjust[#1`#2`#3]{%
\yoff=10
\vertadjust[#1`#2`{#3}]%
\advance \yext by\tempcounta
\advance \yext by 10
}
\def\botadjust[#1`#2`#3]{%
\vertadjust[#1`#2`{#3}]%
\advance \yext by\tempcounta
\advance \yoff by-\tempcounta
}
\def\leftadjust[#1`#2`#3]{%
\xoff=0
\horadjust[#1`#2`{#3}]%
\advance \xext by\tempcounta
\advance \xoff by-\tempcounta
}
\def\rightadjust[#1`#2`#3]{%
\horadjust[#1`#2`{#3}]%
\advance \xext by\tempcounta
}
\def\rightsladjust[#1`#2`#3]{%
\sladjust[#1`#2`{#3}]{\width}%
\advance \xext by\tempcounta
}
\def\leftsladjust[#1`#2`#3]{%
\xoff=0
\sladjust[#1`#2`{#3}]{\width}%
\advance \xext by\tempcounta
\advance \xoff by-\tempcounta
}
\def\adjust[#1`#2;#3`#4;#5`#6;#7`#8]{%
\topadjust[#1``{#2}]
\leftadjust[#3``{#4}]
\rightadjust[#5``{#6}]
\botadjust[#7``{#8}]}
 
\def\putsquarep<#1>(#2)[#3;#4`#5`#6`#7]{{%
\setsqparms[#1]%
\setpos(#2)%
\settokens[#3]%
\puthmorphism(\xpos,\ypos)[\tokenc`\tokend`{#7}]{\width}{\arrowtyped}b%
\advance\ypos by \height
\puthmorphism(\xpos,\ypos)[\tokena`\tokenb`{#4}]{\width}{\arrowtypea}a%
\putvmorphism(\xpos,\ypos)[``{#5}]{\height}{\arrowtypeb}l%
\advance\xpos by \width
\putvmorphism(\xpos,\ypos)[``{#6}]{\height}{\arrowtypec}r%
}}
 
\def\putsquare{\@ifnextchar <{\putsquarep}{\putsquarep%
   <\arrowtypea`\arrowtypeb`\arrowtypec`\arrowtyped;\width`\height>}}
\def\square{\@ifnextchar< {\squarep}{\squarep
   <\arrowtypea`\arrowtypeb`\arrowtypec`\arrowtyped;\width`\height>}}
\def\squarep<#1>[#2`#3`#4`#5;#6`#7`#8`#9]{{
\setsqparms[#1]
\xext=\width                                          
\yext=\height                                         
\topadjust[#2`#3`{#6}]
\botadjust[#4`#5`{#9}]
\leftadjust[#2`#4`{#7}]
\rightadjust[#3`#5`{#8}]
\begin{picture}(\xext,\yext)(\xoff,\yoff)
\putsquarep<\arrowtypea`\arrowtypeb`\arrowtypec`\arrowtyped;\width`\height>%
(0,0)[#2`#3`#4`#5;#6`#7`#8`{#9}]%
\end{picture}%
}}
 
\def\putptrianglep<#1>(#2,#3)[#4`#5`#6;#7`#8`#9]{{%
\settriparms[#1]%
\xpos=#2 \ypos=#3
\advance\ypos by \height
\puthmorphism(\xpos,\ypos)[#4`#5`{#7}]{\height}{\arrowtypea}a%
\putvmorphism(\xpos,\ypos)[`#6`{#8}]{\height}{\arrowtypeb}l%
\advance\xpos by\height
\putmorphism(\xpos,\ypos)(-1,-1)[``{#9}]{\height}{\arrowtypec}r%
}}
 
\def\putptriangle{\@ifnextchar <{\putptrianglep}{\putptrianglep
   <\arrowtypea`\arrowtypeb`\arrowtypec;\height>}}
\def\ptriangle{\@ifnextchar <{\ptrianglep}{\ptrianglep
   <\arrowtypea`\arrowtypeb`\arrowtypec;\height>}}
 
\def\ptrianglep<#1>[#2`#3`#4;#5`#6`#7]{{
\settriparms[#1]%
\width=\height                         
\xext=\width                           
\yext=\width                           
\topadjust[#2`#3`{#5}]
\botadjust[#3``]
\leftadjust[#2`#4`{#6}]
\rightsladjust[#3`#4`{#7}]
\begin{picture}(\xext,\yext)(\xoff,\yoff)
\putptrianglep<\arrowtypea`\arrowtypeb`\arrowtypec;\height>%
(0,0)[#2`#3`#4;#5`#6`{#7}]%
\end{picture}%
}}
 
\def\putqtrianglep<#1>(#2,#3)[#4`#5`#6;#7`#8`#9]{{%
\settriparms[#1]%
\xpos=#2 \ypos=#3
\advance\ypos by\height
\puthmorphism(\xpos,\ypos)[#4`#5`{#7}]{\height}{\arrowtypea}a%
\putmorphism(\xpos,\ypos)(1,-1)[``{#8}]{\height}{\arrowtypeb}l%
\advance\xpos by\height
\putvmorphism(\xpos,\ypos)[`#6`{#9}]{\height}{\arrowtypec}r%
}}
 
\def\putqtriangle{\@ifnextchar <{\putqtrianglep}{\putqtrianglep
   <\arrowtypea`\arrowtypeb`\arrowtypec;\height>}}
\def\qtriangle{\@ifnextchar <{\qtrianglep}{\qtrianglep
   <\arrowtypea`\arrowtypeb`\arrowtypec;\height>}}
 
\def\qtrianglep<#1>[#2`#3`#4;#5`#6`#7]{{
\settriparms[#1]
\width=\height                         
\xext=\width                           
\yext=\height                          
\topadjust[#2`#3`{#5}]
\botadjust[#4``]
\leftsladjust[#2`#4`{#6}]
\rightadjust[#3`#4`{#7}]
\begin{picture}(\xext,\yext)(\xoff,\yoff)
\putqtrianglep<\arrowtypea`\arrowtypeb`\arrowtypec;\height>%
(0,0)[#2`#3`#4;#5`#6`{#7}]%
\end{picture}%
}}
 
\def\putdtrianglep<#1>(#2,#3)[#4`#5`#6;#7`#8`#9]{{%
\settriparms[#1]%
\xpos=#2 \ypos=#3
\puthmorphism(\xpos,\ypos)[#5`#6`{#9}]{\height}{\arrowtypec}b%
\advance\xpos by \height \advance\ypos by\height
\putmorphism(\xpos,\ypos)(-1,-1)[``{#7}]{\height}{\arrowtypea}l%
\putvmorphism(\xpos,\ypos)[#4``{#8}]{\height}{\arrowtypeb}r%
}}
 
\def\putdtriangle{\@ifnextchar <{\putdtrianglep}{\putdtrianglep
   <\arrowtypea`\arrowtypeb`\arrowtypec;\height>}}
\def\dtriangle{\@ifnextchar <{\dtrianglep}{\dtrianglep
   <\arrowtypea`\arrowtypeb`\arrowtypec;\height>}}
 
\def\dtrianglep<#1>[#2`#3`#4;#5`#6`#7]{{
\settriparms[#1]
\width=\height                         
\xext=\width                           
\yext=\height                          
\topadjust[#2``]
\botadjust[#3`#4`{#7}]
\leftsladjust[#3`#2`{#5}]
\rightadjust[#2`#4`{#6}]
\begin{picture}(\xext,\yext)(\xoff,\yoff)
\putdtrianglep<\arrowtypea`\arrowtypeb`\arrowtypec;\height>%
(0,0)[#2`#3`#4;#5`#6`{#7}]%
\end{picture}%
}}
 
\def\putbtrianglep<#1>(#2,#3)[#4`#5`#6;#7`#8`#9]{{%
\settriparms[#1]%
\xpos=#2 \ypos=#3
\puthmorphism(\xpos,\ypos)[#5`#6`{#9}]{\height}{\arrowtypec}b%
\advance\ypos by\height
\putmorphism(\xpos,\ypos)(1,-1)[``{#8}]{\height}{\arrowtypeb}r%
\putvmorphism(\xpos,\ypos)[#4``{#7}]{\height}{\arrowtypea}l%
}}
 
\def\putbtriangle{\@ifnextchar <{\putbtrianglep}{\putbtrianglep
   <\arrowtypea`\arrowtypeb`\arrowtypec;\height>}}
\def\btriangle{\@ifnextchar <{\btrianglep}{\btrianglep
   <\arrowtypea`\arrowtypeb`\arrowtypec;\height>}}
 
\def\btrianglep<#1>[#2`#3`#4;#5`#6`#7]{{
\settriparms[#1]
\width=\height                         
\xext=\width                           
\yext=\height                          
\topadjust[#2``]
\botadjust[#3`#4`{#7}]
\leftadjust[#2`#3`{#5}]
\rightsladjust[#4`#2`{#6}]
\begin{picture}(\xext,\yext)(\xoff,\yoff)
\putbtrianglep<\arrowtypea`\arrowtypeb`\arrowtypec;\height>%
(0,0)[#2`#3`#4;#5`#6`{#7}]%
\end{picture}%
}}
 
\def\putAtrianglep<#1>(#2,#3)[#4`#5`#6;#7`#8`#9]{{%
\settriparms[#1]%
\xpos=#2 \ypos=#3
{\multiply \height by2
\puthmorphism(\xpos,\ypos)[#5`#6`{#9}]{\height}{\arrowtypec}b}%
\advance\xpos by\height \advance\ypos by\height
\putmorphism(\xpos,\ypos)(-1,-1)[#4``{#7}]{\height}{\arrowtypea}l%
\putmorphism(\xpos,\ypos)(1,-1)[``{#8}]{\height}{\arrowtypeb}r%
}}
 
\def\putAtriangle{\@ifnextchar <{\putAtrianglep}{\putAtrianglep
   <\arrowtypea`\arrowtypeb`\arrowtypec;\height>}}
\def\Atriangle{\@ifnextchar <{\Atrianglep}{\Atrianglep
   <\arrowtypea`\arrowtypeb`\arrowtypec;\height>}}
 
\def\Atrianglep<#1>[#2`#3`#4;#5`#6`#7]{{
\settriparms[#1]
\width=\height                         
\xext=\width                           
\yext=\height                          
\topadjust[#2``]
\botadjust[#3`#4`{#7}]
\multiply \xext by2 
\leftsladjust[#3`#2`{#5}]
\rightsladjust[#4`#2`{#6}]
\begin{picture}(\xext,\yext)(\xoff,\yoff)%
\putAtrianglep<\arrowtypea`\arrowtypeb`\arrowtypec;\height>%
(0,0)[#2`#3`#4;#5`#6`{#7}]%
\end{picture}%
}}
 
\def\putAtrianglepairp<#1>(#2)[#3;#4`#5`#6`#7`#8]{{
\settripairparms[#1]%
\setpos(#2)%
\settokens[#3]%
\puthmorphism(\xpos,\ypos)[\tokenb`\tokenc`{#7}]{\height}{\arrowtyped}b%
\advance\xpos by\height
\advance\ypos by\height
\putmorphism(\xpos,\ypos)(-1,-1)[\tokena``{#4}]{\height}{\arrowtypea}l%
\putvmorphism(\xpos,\ypos)[``{#5}]{\height}{\arrowtypeb}m%
\putmorphism(\xpos,\ypos)(1,-1)[``{#6}]{\height}{\arrowtypec}r%
}}
 
\def\putAtrianglepair{\@ifnextchar <{\putAtrianglepairp}{\putAtrianglepairp%
   <\arrowtypea`\arrowtypeb`\arrowtypec`\arrowtyped`\arrowtypee;\height>}}
\def\Atrianglepair{\@ifnextchar <{\Atrianglepairp}{\Atrianglepairp%
   <\arrowtypea`\arrowtypeb`\arrowtypec`\arrowtyped`\arrowtypee;\height>}}
 
\def\Atrianglepairp<#1>[#2;#3`#4`#5`#6`#7]{{%
\settripairparms[#1]%
\settokens[#2]%
\width=\height
\xext=\width
\yext=\height
\topadjust[\tokena``]%
\vertadjust[\tokenb`\tokenc`{#6}]
\tempcountd=\tempcounta                       
\vertadjust[\tokenc`\tokend`{#7}]
\ifnum\tempcounta<\tempcountd                 
\tempcounta=\tempcountd\fi                    
\advance \yext by\tempcounta                  
\advance \yoff by-\tempcounta                 %
\multiply \xext by2 
\leftsladjust[\tokenb`\tokena`{#3}]
\rightsladjust[\tokend`\tokena`{#5}]%
\begin{picture}(\xext,\yext)(\xoff,\yoff)%
\putAtrianglepairp
<\arrowtypea`\arrowtypeb`\arrowtypec`\arrowtyped`\arrowtypee;\height>%
(0,0)[#2;#3`#4`#5`#6`{#7}]%
\end{picture}%
}}
 
\def\putVtrianglep<#1>(#2,#3)[#4`#5`#6;#7`#8`#9]{{%
\settriparms[#1]%
\xpos=#2 \ypos=#3
\advance\ypos by\height
{\multiply\height by2
\puthmorphism(\xpos,\ypos)[#4`#5`{#7}]{\height}{\arrowtypea}a}%
\putmorphism(\xpos,\ypos)(1,-1)[`#6`{#8}]{\height}{\arrowtypeb}l%
\advance\xpos by\height
\advance\xpos by\height
\putmorphism(\xpos,\ypos)(-1,-1)[``{#9}]{\height}{\arrowtypec}r%
}}

\def\putVtriangle{\@ifnextchar <{\putVtrianglep}{\putVtrianglep
   <\arrowtypea`\arrowtypeb`\arrowtypec;\height>}}
\def\Vtriangle{\@ifnextchar <{\Vtrianglep}{\Vtrianglep
   <\arrowtypea`\arrowtypeb`\arrowtypec;\height>}}
 
\def\Vtrianglep<#1>[#2`#3`#4;#5`#6`#7]{{
\settriparms[#1]
\width=\height                         
\xext=\width                           
\yext=\height                          
\topadjust[#2`#3`{#5}]
\botadjust[#4``]
\multiply \xext by2 
\leftsladjust[#2`#3`{#6}]
\rightsladjust[#3`#4`{#7}]
\begin{picture}(\xext,\yext)(\xoff,\yoff)%
\putVtrianglep<\arrowtypea`\arrowtypeb`\arrowtypec;\height>%
(0,0)[#2`#3`#4;#5`#6`{#7}]%
\end{picture}%
}}
 
\def\putVtrianglepairp<#1>(#2)[#3;#4`#5`#6`#7`#8]{{
\settripairparms[#1]%
\setpos(#2)%
\settokens[#3]%
\advance\ypos by\height
\putmorphism(\xpos,\ypos)(1,-1)[`\tokend`{#6}]{\height}{\arrowtypec}l%
\puthmorphism(\xpos,\ypos)[\tokena`\tokenb`{#4}]{\height}{\arrowtypea}a%
\advance\xpos by\height
\putvmorphism(\xpos,\ypos)[``{#7}]{\height}{\arrowtyped}m%
\advance\xpos by\height
\putmorphism(\xpos,\ypos)(-1,-1)[``{#8}]{\height}{\arrowtypee}r%
}}
 
\def\putVtrianglepair{\@ifnextchar <{\putVtrianglepairp}{\putVtrianglepairp%
    <\arrowtypea`\arrowtypeb`\arrowtypec`\arrowtyped`\arrowtypee;\height>}}
\def\Vtrianglepair{\@ifnextchar <{\Vtrianglepairp}{\Vtrianglepairp%
    <\arrowtypea`\arrowtypeb`\arrowtypec`\arrowtyped`\arrowtypee;\height>}}
 
\def\Vtrianglepairp<#1>[#2;#3`#4`#5`#6`#7]{{%
\settripairparms[#1]%
\settokens[#2]
\xext=\height                  
\width=\height                 
\yext=\height                  
\vertadjust[\tokena`\tokenb`{#4}]
\tempcountd=\tempcounta        
\vertadjust[\tokenb`\tokenc`{#5}]
\ifnum\tempcounta<\tempcountd%
\tempcounta=\tempcountd\fi
\advance \yext by\tempcounta
\botadjust[\tokend``]%
\multiply \xext by2
\leftsladjust[\tokena`\tokend`{#6}]%
\rightsladjust[\tokenc`\tokend`{#7}]%
\begin{picture}(\xext,\yext)(\xoff,\yoff)%
\putVtrianglepairp
<\arrowtypea`\arrowtypeb`\arrowtypec`\arrowtyped`\arrowtypee;\height>%
(0,0)[#2;#3`#4`#5`#6`{#7}]%
\end{picture}%
}}

\def\putCtrianglep<#1>(#2,#3)[#4`#5`#6;#7`#8`#9]{{%
\settriparms[#1]%
\xpos=#2 \ypos=#3
\advance\ypos by\height
\putmorphism(\xpos,\ypos)(1,-1)[``{#9}]{\height}{\arrowtypec}l%
\advance\xpos by\height
\advance\ypos by\height
\putmorphism(\xpos,\ypos)(-1,-1)[#4`#5`{#7}]{\height}{\arrowtypea}l%
{\multiply\height by 2
\putvmorphism(\xpos,\ypos)[`#6`{#8}]{\height}{\arrowtypeb}r}%
}}
 
\def\putCtriangle{\@ifnextchar <{\putCtrianglep}{\putCtrianglep
    <\arrowtypea`\arrowtypeb`\arrowtypec;\height>}}
\def\Ctriangle{\@ifnextchar <{\Ctrianglep}{\Ctrianglep
    <\arrowtypea`\arrowtypeb`\arrowtypec;\height>}}
 
\def\Ctrianglep<#1>[#2`#3`#4;#5`#6`#7]{{
\settriparms[#1]
\width=\height                          
\xext=\width                            
\yext=\height                           
\multiply \yext by2 
\topadjust[#2``]
\botadjust[#4``]
\sladjust[#3`#2`{#5}]{\width}
\tempcountd=\tempcounta                 
\sladjust[#3`#4`{#7}]{\width}
\ifnum \tempcounta<\tempcountd          
\tempcounta=\tempcountd\fi              
\advance \xext by\tempcounta            
\advance \xoff by-\tempcounta           %
\rightadjust[#2`#4`{#6}]
\begin{picture}(\xext,\yext)(\xoff,\yoff)%
\putCtrianglep<\arrowtypea`\arrowtypeb`\arrowtypec;\height>%
(0,0)[#2`#3`#4;#5`#6`{#7}]%
\end{picture}%
}}
 
\def\putDtrianglep<#1>(#2,#3)[#4`#5`#6;#7`#8`#9]{{%
\settriparms[#1]%
\xpos=#2 \ypos=#3
\advance\xpos by\height \advance\ypos by\height
\putmorphism(\xpos,\ypos)(-1,-1)[``{#9}]{\height}{\arrowtypec}r%
\advance\xpos by-\height \advance\ypos by\height
\putmorphism(\xpos,\ypos)(1,-1)[`#5`{#8}]{\height}{\arrowtypeb}r%
{\multiply\height by 2
\putvmorphism(\xpos,\ypos)[#4`#6`{#7}]{\height}{\arrowtypea}l}%
}}
 
\def\putDtriangle{\@ifnextchar <{\putDtrianglep}{\putDtrianglep
    <\arrowtypea`\arrowtypeb`\arrowtypec;\height>}}
\def\Dtriangle{\@ifnextchar <{\Dtrianglep}{\Dtrianglep
   <\arrowtypea`\arrowtypeb`\arrowtypec;\height>}}
 
\def\Dtrianglep<#1>[#2`#3`#4;#5`#6`#7]{{
\settriparms[#1]
\width=\height                         
\xext=\height                          
\yext=\height                          
\multiply \yext by2 
\topadjust[#2``]
\botadjust[#4``]
\leftadjust[#2`#4`{#5}]
\sladjust[#3`#2`{#5}]{\height}
\tempcountd=\tempcountd                
\sladjust[#3`#4`{#7}]{\height}
\ifnum \tempcounta<\tempcountd         
\tempcounta=\tempcountd\fi             
\advance \xext by\tempcounta           %
\begin{picture}(\xext,\yext)(\xoff,\yoff)
\putDtrianglep<\arrowtypea`\arrowtypeb`\arrowtypec;\height>%
(0,0)[#2`#3`#4;#5`#6`{#7}]%
\end{picture}%
}}
 
\def\setrecparms[#1`#2]{\width=#1 \height=#2}%
\def\recursep<#1`#2>[#3;#4`#5`#6`#7`#8]{{%
\width=#1 \height=#2
\settokens[#3]
\settowidth{\tempdimen}{$\tokena$}
\ifdim\tempdimen=0pt
  \savebox{\tempboxa}{\hbox{$\tokenb$}}%
  \savebox{\tempboxb}{\hbox{$\tokend$}}%
  \savebox{\tempboxc}{\hbox{$#6$}}%
\else
  \savebox{\tempboxa}{\hbox{$\hbox{$\tokena$}\times\hbox{$\tokenb$}$}}%
  \savebox{\tempboxb}{\hbox{$\hbox{$\tokena$}\times\hbox{$\tokend$}$}}%
  \savebox{\tempboxc}{\hbox{$\hbox{$\tokena$}\times\hbox{$#6$}$}}%
\fi
\ypos=\height
\divide\ypos by 2
\xpos=\ypos
\advance\xpos by \width
\xext=\xpos \yext=\height
\topadjust[#3`\usebox{\tempboxa}`{#4}]%
\botadjust[#5`\usebox{\tempboxb}`{#8}]%
\sladjust[\tokenc`\tokenb`{#5}]{\ypos}%
\tempcountd=\tempcounta
\sladjust[\tokenc`\tokend`{#5}]{\ypos}%
\ifnum \tempcounta<\tempcountd
\tempcounta=\tempcountd\fi
\advance \xext by\tempcounta
\advance \xoff by-\tempcounta
\rightadjust[\usebox{\tempboxa}`\usebox{\tempboxb}`\usebox{\tempboxc}]%
\bfig
\putCtrianglep<-1`1`1;\ypos>(0,0)[`\tokenc`;#5`#6`{#7}]%
\puthmorphism(\ypos,0)[\tokend`\usebox{\tempboxb}`{#8}]{\width}{-1}b%
\puthmorphism(\ypos,\height)[\tokenb`\usebox{\tempboxa}`{#4}]{\width}{-1}a%
\advance\ypos by \width
\putvmorphism(\ypos,\height)[``\usebox{\tempboxc}]{\height}1r%
\efig
}}
 
\def\recurse{\@ifnextchar <{\recursep}{\recursep<\width`\height>}}
 
\def\puttwohmorphisms(#1,#2)[#3`#4;#5`#6]#7#8#9{{%
\puthmorphism(#1,#2)[#3`#4`]{#7}0a
\ypos=#2
\advance\ypos by 20
\puthmorphism(#1,\ypos)[\phantom{#3}`\phantom{#4}`#5]{#7}{#8}a
\advance\ypos by -40
\puthmorphism(#1,\ypos)[\phantom{#3}`\phantom{#4}`#6]{#7}{#9}b
}}
 
\def\puttwovmorphisms(#1,#2)[#3`#4;#5`#6]#7#8#9{{%
\putvmorphism(#1,#2)[#3`#4`]{#7}0a
\xpos=#1
\advance\xpos by -20
\putvmorphism(\xpos,#2)[\phantom{#3}`\phantom{#4}`#5]{#7}{#8}l
\advance\xpos by 40
\putvmorphism(\xpos,#2)[\phantom{#3}`\phantom{#4}`#6]{#7}{#9}r
}}
 
\def\puthcoequalizer(#1)[#2`#3`#4;#5`#6`#7]#8#9{{%
\setpos(#1)%
\puttwohmorphisms(\xpos,\ypos)[#2`#3;#5`#6]{#8}11%
\advance\xpos by #8
\puthmorphism(\xpos,\ypos)[\phantom{#3}`#4`#7]{#8}1{#9}
}}
 
\def\putvcoequalizer(#1)[#2`#3`#4;#5`#6`#7]#8#9{{%
\setpos(#1)%
\puttwovmorphisms(\xpos,\ypos)[#2`#3;#5`#6]{#8}11%
\advance\ypos by -#8
\putvmorphism(\xpos,\ypos)[\phantom{#3}`#4`#7]{#8}1{#9}
}}
 
\def\putthreehmorphisms(#1)[#2`#3;#4`#5`#6]#7(#8)#9{{%
\setpos(#1) \settypes(#8)
\if a#9 %
     \vertsize{\tempcounta}{#5}%
     \vertsize{\tempcountb}{#6}%
     \ifnum \tempcounta<\tempcountb \tempcounta=\tempcountb \fi
\else
     \vertsize{\tempcounta}{#4}%
     \vertsize{\tempcountb}{#5}%
     \ifnum \tempcounta<\tempcountb \tempcounta=\tempcountb \fi
\fi
\advance \tempcounta by 60
\puthmorphism(\xpos,\ypos)[#2`#3`#5]{#7}{\arrowtypeb}{#9}
\advance\ypos by \tempcounta
\puthmorphism(\xpos,\ypos)[\phantom{#2}`\phantom{#3}`#4]{#7}{\arrowtypea}{#9}
\advance\ypos by -\tempcounta \advance\ypos by -\tempcounta
\puthmorphism(\xpos,\ypos)[\phantom{#2}`\phantom{#3}`#6]{#7}{\arrowtypec}{#9}
}}
 
\def\putarc(#1,#2)[#3`#4`#5]#6#7#8{{%
\xpos #1
\ypos #2
\width #6
\arrowlength #6
\putbox(\xpos,\ypos){#3\vphantom{#4}}%
{\advance \xpos by\arrowlength
\putbox(\xpos,\ypos){\vphantom{#3}#4}}%
\horsize{\tempcounta}{#3}%
\horsize{\tempcountb}{#4}%
\divide \tempcounta by2
\divide \tempcountb by2
\advance \tempcounta by30
\advance \tempcountb by30
\advance \xpos by\tempcounta
\advance \arrowlength by-\tempcounta
\advance \arrowlength by-\tempcountb
\halflength=\arrowlength \divide\halflength by 2
\divide\arrowlength by 5
\put(\xpos,\ypos){\bezier{\arrowlength}(0,0)(50,50)(\halflength,50)}
\ifnum #7=-1 \put(\xpos,\ypos){\vector(-3,-2)0} \fi
\advance\xpos by \halflength
\put(\xpos,\ypos){\xpos=\halflength \advance\xpos by -50
   \bezier{\arrowlength}(0,50)(\xpos,50)(\halflength,0)}
\ifnum #7=1 {\advance \xpos by
   \halflength \put(\xpos,\ypos){\vector(3,-2)0}} \fi
\advance\ypos by 50
\vertsize{\tempcounta}{#5}%
\divide\tempcounta by2
\advance \tempcounta by20
\if a#8 %
   \advance \ypos by\tempcounta
   \putbox(\xpos,\ypos){#5}%
\else
   \advance \ypos by-\tempcounta
   \putbox(\xpos,\ypos){#5}%
\fi
}}
 
\makeatother

\input{Macro}
\input{macrosTikz}

\def\lastname{S. Ghilardi and L. Santocanale}

\ifams
\author[Ghilardi]{Silvio Ghilardi}
\address{Silvio Ghilardi\\
 Dipartimento di Matematica, Universit\`a degli Studi di
  Milano}
\email{silvio.ghilardi@unimi.it}
\author[Santocanale]{Luigi Santocanale}
\address{Luigi Santocanale\\
LIS, CNRS UMR 7020, Aix-Marseille Universit\'e}
\email{luigi.santocanale@lis-lab.fr}
\fi

\ifmscs
\author[\lastname]{Silvio Ghilardi$^{1}$ and Luigi Santocanale$^{2}$\\
  $^1$ Dipartimento di Matematica, Universit\`a degli Studi di
  Milano 
  \addressbreak
  $^2$ LIS, CNRS UMR 7020, Aix-Marseille Universit\'e
}
\fi

\title[Free Heyting Algebra Endomorphisms ...]{Free Heyting Algebra Endomorphisms: \\ Ruitenburg's Theorem and Beyond}

\begin{document}

\maketitle

\begin{abstract}
  \RT says that every endomorphism $f$ of a finitely generated free
  Heyting algebra is ultimately periodic if $f$ fixes all the
  generators but one. More precisely, there is $N\geq 0$ such that
  $f^{N+2}= f^N$, thus the period equals 2.  We give a semantic proof
  of this theorem, using duality techniques and bounded bisimulation
  ranks.  By the same techniques, we tackle investigation of arbitrary
  endomorphisms between free algebras. We show that they are not, in
  general, ultimately periodic. Yet, when they are (e.g. in the case
  of locally finite subvarieties), the period can be explicitly
  bounded as function of the cardinality of the set of generators.


\end{abstract}

\section{Introduction}\label{sec:intro}

Unification theory investigates the behavior of substitutions from a
syntactic point of view: substitutions are in fact key ingredients in
various algorithms commonly used in computational logic. Taking an
algebraic point of view, substitutions can be seen as finitely
generated free algebra homomorphisms:
in fact, such a homomorphism 
$$\mu: \cF(x_1, \dots, x_n)\lora \cF(y_1, \dots, y_m)$$
is uniquely determined by an $n$-tuple of terms 
$$t_1(y_1, \dots, y_m), \dots, t_n(y_1, \dots, y_m)$$ 
and acts by
associating with any term $u(x_1, \dots, x_n)$ the term
$$u(t_1/x_1, \dots, t_n/x_n)$$ 
obtained by substitution. If free algebras are intended not as
`absolutely free algebras', but as `free algebras in an equational
class $E$', the same correspondence between homomorphisms and
substitutions works, provided terms are intended as equivalence
classes of terms modulo $E$ and substitutions themselves are taken
`modulo $E$'.

The above correspondence between free
  algebra homomorphisms and substitutions is the starting point for
the algebraic approaches to $E$-uni\-fi\-ca\-tion theory,
like for instance~\cite{goguen,unpro}, where structural information about
homomorphisms of finitely generated (and also finitely presented)
algebras is widely exploited. In this paper, we want to draw the
attention on a surprising behavior that such homomorphisms can have in
some algebraic logic contexts. Such behavior is unexpectedly similar
to that of functions between finite sets.

To explain what we have in mind, let us recall that an infinite sequence
$$
a_1, a_2, \dots, a_i, \dots
$$
\emph{ultimately periodic} if there are $N$ and $k$
such that for all $s_1, s_2\geq N$, we have that
$s_1\equiv s_2 \mod k$ implies $a_{s_1}= a_{s_2}$.  If $(N, k)$ is the
smallest (in the lexicographic sense) pair for which this happens,
then $N$ and $k$ are, respectively, the \emph{index} and the
  \emph{period} of the ultimately periodic sequence
$\set{a_i}_i$. Thus, for instance, an ultimately periodic sequence
with index $N$ and period $2$ looks as follows
$$
a_1, \dots, a_N, a_{N+1}, a_N, a_{N+1}, \dots
$$
A typical example of an ultimately periodic sequence is the sequence
of the iterations $\set{f^i}_i$ of an endo-function $f$ of a finite
set.  Whenever infinitary data are involved, ultimate periodicity comes
often as a surprise.

  \RT is in fact a surprising result stating the following: take a
  formula $A(x, \uy)$ of intuitionistic propositional calculus $(IPC)$
  (by the notation $A(x, \uy)$ we mean that the only propositional
  letters occurring in $A$ are among $x, \uy$ - with $\uy$ being, say,
  the tuple $y_1, \dots, y_n$) and consider the sequence
  $\set{A^i(x,\uy)}_{i\geq 1}$ so defined:
\begin{equation}\label{eq:formulasequence}
A^1\eqdef A, ~~\dots,~~ A^{i+1}\eqdef A(A^i/x, \uy)
\end{equation}
where the
  slash means substitution; then,
  \emph{taking equivalence
  classes under provable bi-implication in $(IPC)$, the sequence
  $\set{[A^i(x,\uy)]}_{i\geq 1}$ is ultimately periodic with period
  2}.
  The latter means that there is $N$ such that
\begin{equation}\label{eq:Ruitenburgtheorem}
\vdash_{IPC} A^{N+2}\leftrightarrow  A^N~~~.
\end{equation}

An interesting consequence of this result is that \emph {least (and
  greatest) fixpoints of monotonic formulae are definable in
  $(IPC)$}~\cite{Mardaev1993,Mardaev07,fossacs}: this is because the
sequence~\eqref{eq:formulasequence} becomes increasing when evaluated
on $\bot/x$ (if $A$ is monotonic in $x$), so that the period is
decreased to 1.  Thus the index of the sequence becomes a finite upper
bound for the fixpoint approximations convergence: in fact we have,
$\vdash_{IPC} A^N(\bot/x) \to A^{N+1}(\bot/x)$ and
$\vdash_{IPC} A^{N+1}(\bot/x) \to A^{N+2}(\bot/x)$ by the monotonicity
of $A$, yielding
$\vdash_{IPC} A^N(\bot/x) \leftrightarrow A^{N+1}(\bot/x)$
by~\eqref{eq:Ruitenburgtheorem}.

\RT was shown in~\cite{Ruitenburg84} via a, rather involved, purely
syntactic proof. The proof has been recently formalized inside the
proof assistant \textsc{coq} by T. Litak, see 
\url{https://git8.cs.fau.de/redmine/projects/ruitenburg1984}~~.
In
this paper we supply a semantic proof, using duality and bounded
bisimulation machinery.

\emph{Bounded bisimulations} are a standard tool in non classical
logics~\cite{fine} which is used in order to characterize
satisfiability of bounded depth formulae and hence definable classes
of models: examples of the use of bounded bisimulations include for
instance~\cite{shavrukov,GhilardiZawadowski2011,visser,um}.

\emph{Duality} has a long tradition in algebraic logic, see
e.g.~\cite{Esa74} for the case of \Ha{s}. Indeed, many phenomena look
more transparent whenever they are analyzed in the dual
categories. This especially happens when dualities can convert
coproducts and colimits constructions into more familiar `honest'
products and limits constructions.  The duality we use to tackle \RT,
firstly described in~\cite{GhilardiZawadowski97}, see
also~\cite{GhilardiZawadowski2011}, realizes this conversion.  It has
a mixed geometric/combinatorial nature.
In fact, the geometric environment shows \emph{how to find} relevant
mathematical structures (products, equalizers, images,...) using their
standard definitions in sheaves and presheaves; on the other hand, the
combinatorial aspects show that such constructions \emph{are
  definable}, thus meaningful from the logical side. In this sense,
notice that we work with finitely presented algebras, and our
combinatorial ingredients (Ehrenfeucht-Fraiss\'e games, etc.) replace
the topological ingredients which are common in the algebraic logic
literature (working with arbitrary algebras instead).  Duality,
although not always in an explicitly mentioned form, is also at the
heart of the finitarity results for $E$-unification theory
in~\cite{G99,um,unapal}.

The paper is organized as follows. In Section~\ref{sec:classical} we
show how to formulate \RT in algebraic terms and how to prove it via
duality in the easy case of classical logic (where index is always
1). This Section supplies the methodology we shall follow in the whole
paper. 
We introduce in Section~\ref{sec:duality} the required duality
ingredients for finitely presented \Ha{s}, leading to the statement of
the duality Theorem. The full proof of this theorem appears in the
following Section~\ref{sec:dualityproof}.  We show then, in
Section~\ref{sec:finite}, how to extend the basic argument of
Section~\ref{sec:classical} to finite Kripke models of \IL. 
This extension does not directly give \RT, because it supplies a bound
for the indexes of our sequences which is dependent on the poset a
given model is based on. Using the ranks machinery introduced in
Section~\ref{sec:ranks}, this bound is made uniform in
Section~\ref{sec:main}, thus finally reaching our first goal.  Having
established \RT, 
we wonder how general this ultimately periodic behavior is among the
finitely generated free \Ha endomorphisms and, in
Section~\ref{sec:endomorphisms}, we supply a counterexample showing
that this behavior fails whenever at least two free generators are
moved by the endomorphism. In the final Section~\ref{sec:bounds}, we
prove that, whenever an endomorphism is ultimately periodic, its
period can be bound as a function of the number of the free generators
only. This observation is used to provide bounds of periods of free
algebra endomorphisms in locally finite
varieties of \Ha{s}.
We present concluding remarks and some open problems in the last
Section.

Most of the material of this paper was presented at the conference
AiML 18, see the reference~\cite{aiml18}; the content of the last two
Sections as well as a strengthening of the duality theorem of
\cite{GhilardiZawadowski97}, however, are novel.

\section{The Case of Classical Logic}\label{sec:classical}

We explain our methodology in the much easier case of classical
logic. In classical propositional calculus ($CPC$), \RT holds with
index 1 and period 2, namely given a formula $A(x,\uy)$, we prove that
\begin{equation}
  \label{eq:class}
 \vdash_{CPC} A^3 \leftrightarrow A 
\end{equation}
holds (here $A^3$ is defined like in~\eqref{eq:formulasequence}). 
\subsection{The algebraic reformulation}

First, we transform the above statement \eqref{eq:class}
into an algebraic statement concerning
free Boolean algebras. We let $\cF_B(\uz)$ be the free Boolean algebra
over the finite set $\uz$.  Recall that $\cF_B(\uz)$ is the
Lindenbaum-Tarski algebra of classical propositional calculus
restricted to a language having just the $\uz$ as propositional
variables.

Similarly, morphisms $\mu: \cF_B(x_1, \dots, x_n)\lora \cF_B(\uz)$
bijectively correspond to $n$-tuples of equivalence classes of
formulae $A_1(\uz), \dots, A_n(\uz)$ in $\cF_B(\uz)$: the map $\mu$
corresponding to the tuple $A_1(\uz), \dots, A_n(\uz)$ associates with
the equivalence class of $B(x_1, \dots, x_n)$ in
$\cF_B(x_1, \dots, x_n)$ the equivalence class of
$B(A_1/x_1, \dots, A_n/x_n)$ in $\cF_B(\uz)$.

Composition is substitution, in the sense that if
$\mu: \cF_B(x_1, \dots, x_n)\lora \cF_B(\uz)$ is induced, as above, by
$A_1(\uz), \dots, A_n(\uz)$ and if
$\nu:\cF_B(y_1, \dots, y_m) \lora \cF_B(x_1, \dots, x_n)$ is induced
by $C_1(x_1, \dots, x_n), \dots, C_m(x_1, \dots, x_n)$, then the map
$\mu\circ \nu: \cF_B(y_1, \dots, y_m)\lora \cF_B(\uz)$ is induced by
the $m$-tuple
\ifams
of formulas
\fi
$C_1(A_1/x_1, \dots, A_n/x_n), \dots, C_m(A_1/x_1, \dots, A_n/x_n)$.

How to translate the statement~\eqref{eq:class} in this setting? Let
$\uy$ be $y_1, \dots, y_n$; we can consider the map
$\mu_A:\cF_B(x,y_1, \dots, y_n)\lora \cF_B(x,y_1, \dots, y_n)$ induced
by the $n+1$-tuple of formulae $A, y_1, \dots, y_n$; then, taking in
mind that in Lindenbaum algebras identity is modulo provable
equivalence, the statement~\eqref{eq:class} is equivalent to
\begin{equation}\label{eq:class1}
 \mu_A^3 = \mu_A~~. 
\end{equation}
This raises the question: which endomorphisms of $\cF_B(x,\uy)$ are of the kind $\mu_A$ for some $A(x, \uy)$? The answer is simple: consider the `inclusion'
map $\iota$ of  $\cF_B(\uy)$ into $\cF_B(x,\uy)$ (this is the map induced by the $n$-tuple $y_1, \dots, y_n$): the maps 
$\mu:\cF_B(x,\uy)\lora\cF_B(x,\uy) $
that are of the kind $\mu_A$ are precisely the maps $\mu$ such that $\mu\circ \iota= \iota$, i.e. those for which the triangle
\begin{center}
\resetparms
\settriparms[1`1`1;400] \Atriangle[\cF_B(\uy)`\cF_B(x,\uy) 
`\cF_B(x,\uy) ;\iota`\iota`\mu ]
\end{center}
\noindent
commutes. 

It is worth making a little step further: since the free algebra
functor preserves coproducts, we have that $\cF_B(x,\uy)$ is the
coproduct of $\cF_B(\uy)$ with $\cF_B(x)$ - the latter being the free
algebra on one generator. In general, let us denote by $\cA[x]$ the
coproduct of the Boolean algebra $\cA$ with the free algebra on one
generator (let us call $\cA[x]$ the \emph{algebra of polynomials} over
$\cA$). 

 Recall that an algebra is
\emph{finitely presented} if it is isomorphic to the
quotient of a finitely generated free algebra by a finitely generated
congruence. For Boolean algebras, being
  `finitely presented' is equivalent to being `finite'. Yet, we should
  keep mind in the following Sections that this equivalence fails for
  \Ha{s}---so the two notions are in general distinct. 
A slight generalization of statement~\eqref{eq:class1} now reads as
follows:
\begin{itemize}
\item let $\cA$ be a finitely presented Boolean
  algebra 
   and let the
  map $\mu: \cA[x]\lora \cA[x]$ commute with the coproduct injection
  $\iota: \cA \lora \cA[x]$
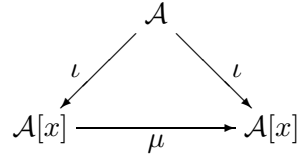
\begin{center}
\resetparms
\settriparms[1`1`1;400] \Atriangle[\cA`\cA[x] 
`\cA[x] ;\iota`\iota`\mu ]
\end{center}
\noindent
Then we have 
\begin{equation}\label{eq:class2}
\mu^3=\mu~~. 
\end{equation}
\end{itemize}

\subsection{Duality}

The gain we achieved with statement~\eqref{eq:class2} is that the latter is a purely categorical statement, so that we can re-interpret it in dual categories.
In fact, a good duality may turn coproducts into products and make our statement easier - if not trivial at all.

Finitely presented Boolean algebras are dual to finite sets; the duality functor maps coproducts into products and the free Boolean algebra on one generator 
to the two-elements set ${\bf 2}=\set{0,1}$ (which, by chance is also a subobject classifier for finite sets). Thus statement~\eqref{eq:class2} now becomes
\begin{itemize}
\item let $T$ be a finite set and let the function
  $f: T\times {\bf{2}}\lora T\times {\bf{2}}$ commute with the product
  projection $\pi_0: T\times {\bf{2}} \lora T$
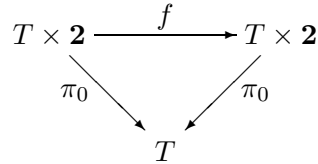
  \begin{center}
    \resetparms \settriparms[1`1`1;400] \Vtriangle[ T\times {\bf{2}}`
    T\times {\bf{2}} `T ;f`\pi_0`\pi_0 ]
\end{center}
\noindent
Then we have 
\begin{equation}\label{eq:class3}
f^3=f~~. 
\end{equation}
\end{itemize}

In this final form, statement~\eqref{eq:class3} is now just a trivial
exercise, which is solved as follows. Notice first that $f$ can be
decomposed as $\langle \pi_0, \chi_S\rangle$ (incidentally, $\chi_S$
is the characteristic function of some $S\subseteq T\times \bf{2}$).
Now, if $f(a,b)=(a,b)$ we trivially have also $f^3(a,b)=f(a,b)$;
suppose then $f(a,b)=(a, b')\neq (a,b)$.  If $f(a,b')=(a,b')$, then
$f^3(a,b)=f(a,b)=(a,b')$, otherwise $f(a,b')=(a,b)$ (there are only
two available values for $b$!) and even in this case
$f^3(a,b)=f(a,b)$.

 Let us illustrate theses cases by thinking of $f$ as an action of the
 monoid of natural numbers 
 on the set $A\times {\bf 2}$, that is, as one-letter deterministic
 automaton:
 
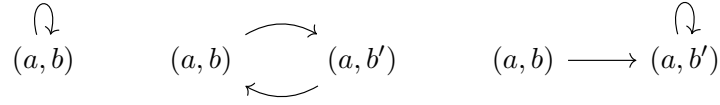
\begin{center}
   \begin{tikzcd}
     (a,b)\ar[loop above]{}{} &
     (a,b) \ar[r, bend left]{}{} & (a,b') \ar[l,, bend left]{}{}&
     (a,b) \ar[r]{}{} & (a,b') \ar[loop above]{}{}
  \end{tikzcd}
 \end{center}
 On each connected component of the automaton,
 the pair index/period is among $(0,1)$, $(0,2)$, $(1,1)$.
 We can compute the global index/period of $f$ by means of
 a $\max/\lcm$ formula: $(1,2) = (\max\set{0,0,1},\lcm\set{1,2})$.

\section{Duality for Heyting Algebras}
\label{sec:duality}

In this Section we supply definitions, notation and statements
from~\cite{GhilardiZawadowski2011} concerning duality for \fp Heyting
algebras.  

A partially ordered set (poset, for short) is a set endowed with a
reflexive, transitive, antisymmetric relation (to be always denoted
with $\leq$).  A poset $P$ is rooted if it has a greatest element,
that we shall denote by $\Root{P}$.  If a finite poset $L$ is fixed,
we call an \emph{$L$-evaluation} or simply an \emph{evaluation} a pair
$\langle P, u\rangle$, where $P$ is a rooted finite poset and
$u:P\ra L$ is an order-preserving map.

Evaluations \emph{restrictions} are introduced as follows.
If $\langle P, u\rangle$ is an $L$-evaluation and if $p \in P$,
then we shall denote by $u_{p}$ the  $L$-evaluation
$\langle \downset p, u \circ i\rangle$, where
$\downset p = \set{p ' \in P \mid p' \leq p}$ and
$i : \,\downset p \subseteq P$ is the inclusion map; briefly,
$u_{p}$ is the restriction of $u$ to the downset generated by $p$.

Evaluations have a strict relationship with finite Kripke models: we
show in detail the connection.  If $\langle L, \leq\rangle$ is
$\langle {\cal P}(\ux), \supseteq\rangle$ (where $\ux=x_1, \dots, x_n$
is a finite list of propositional letters), then an $L$-evaluation
$u:P\ra L$ is called a \emph{Kripke model} for the propositional
intuitionistic language built up from $\ux$.\footnote{
  However, let us
    notice that, according to our convention, a
    $\langle P(\vec{x}),\supseteq\rangle$-evaluation is such that, for
    $p, q\in P$ if $p\leq q$ then $u(p)\supseteq u(q)$; in standard
    logical literature, see e.g. \cite{CZ}, the opposite order on $P$ is
    used, namely an evaluation is such that $u(q)\subseteq u(p)$, for
    $q \leq p$.} 
Given such a Kripke model $u$ and an IPC formula $A(\ux)$, the
\emph{forcing} relation $u\models A$ is inductively defined as
follows:
$$
\begin{aligned}
 &u\models x_i~ &{\rm iff}~ &x_i\in u(\Root{P})
  \\
 &u\not \models \bot &&
 \\
 &u\models A_1\wedge A_2~&{\rm iff}~&(u\models A_1~{\rm and}~u\models A_2)
 \\
 &u\models A_1\vee A_2~&{\rm iff}~&(u\models A_1~{\rm or}~u\models A_2)
 \\
 &u\models A_1\to A_2~&{\rm iff}~&\forall q\leq \Root{P}~(u_q\models A_1~{\Rightarrow}~u_q\models A_2)~~.
\end{aligned}
$$

We define for every $n\in \omega$ and for every pair of
$L$-evaluations $u$ and $v$, the notions of being {\it $n$-equivalent}
(written $u\sim_n v$). 
We also define, for two
  $L$-evaluations $u, v$, the notions of being {\it infinitely
    equivalent} (written $u\sim_{\infty}v$).

Let $u:P\ra L$ and $v:Q\ra L$ be two $L$-evaluations. 
The {\it game} we are interested in has two
  players, Player 1 and Player 2.  
Player 1 can choose either a point in $P$ or a point in $Q$ and Player
2 must answer by choosing a point in the other poset; the only rule of the game is that, if
$\langle p\in P, q\in Q\rangle$ is the last move played so far, then
in the successive move the two players can only choose points
$\langle p', q'\rangle$ such that $p'\leq p$ and $q'\leq q$. If
$\langle p_1, q_1\rangle, \dots, \langle p_i, q_i\rangle, \dots$ are
the points chosen in the game, Player 2 wins iff for every
$i=1, 2, \dots$, we have that $u(p_i)=v(q_i)$. We say that
\begin{itemize}
\item[-]
$u\sim_{\infty} v$ iff {\it Player 2 has a winning strategy} in the above game with infinitely many moves;
\item[-]
$u\sim_n v$ (for $n>0$) iff {\it Player 2 has a winning strategy}
in the above game with $n$ moves, i.e. he has a winning strategy provided
we stipulate that the game terminates after $n$ moves;
\item[-] $u\sim_0 v$ iff $u(\Root{P})=v(\Root{Q})$
  (recall that $\Root{P}, \Root{Q}$ denote the roots
  of $P, Q$).
\end{itemize}

Notice that $u\sim_n v$ always implies $u\sim_0 v$, by the fact that
$L$-evaluations are order-preserving. 
We shall use the notation $[v]_{n}$ for the equivalence class of
  an $L$-valuation $v$ via the equivalence relation $\sim_{n}$.

    The following Proposition provides an elementary recursive
    characterization of the relations $\sim_{n}$, $n \geq 1$. Keeping
    the above definition for $\sim_0$ as base case for recursion, the
    Proposition supplies an alternative recursive definition for these
    relations.  
\begin{proposition}\label{p41.1}  
 Given two $L$-evaluations
$u:P\ra L, v:Q\ra L$, and $n>0$, we have that
$u\sim_{n+1}v$ iff $\forall p\in P\;\exists q\in Q ~(u_p\sim_n
v_q)$  and vice versa.
\end{proposition}

When $L={\cal P}(x_1, \dots, x_n)$, so $L$-evaluations are just
ordinary finite Kripke models over the language built up from the
propositional variables $x_1, \dots, x_n$, the relations $\sim_{n}$
are related to the implicational degree of formulas.
For an IPC formula $A(\ux)$, its implicational degree $d(A)$ is
defined as follows:
\begin{description}
 \item[{\rm (i)}] $d(\bot)=d(x_i)=0$, for $x_i\in \ux$;
 \item[{\rm (ii)}] $d(A_1*A_2)= max[d(A_1), d(A_2)]$, for $*=\wedge, \vee$;
 \item[{\rm (iii)}] $d(A_1\to A_2)= max[d(A_1), d(A_2)]+1$.
\end{description}
One can prove~\cite{visser} that: (1) $u\sim_\infty v$ holds precisely
when ($u\models A \Leftrightarrow v\models A$) holds for all formulae
$A(\ux)$; (2) for all $n$, $u\sim_n v$ holds precisely when
($u\models A \Leftrightarrow v\models A$) holds for all formulae
$A(\ux)$ with $d(A)\leq n$. That is, two evaluations are
$\sim_{\infty}$-equivalent iff they force the same formulas and they
are $\sim_n$-equivalent iff they force the same formulas up to
implicational degree $n$.  Let us remark that, for (1) to be true, it
is essential that our evaluations are defined over \emph{finite}
posets.

The above discussion motivates a sort of  identification of formulae with sets of evaluations closed under restrictions and under $\sim_n$ for some $n$. Thus,
\emph{bounded bisimulations} (this is the way the relations $\sim_n$
are sometimes called) supply the combinatorial ingredients for our
duality; for the picture to be complete, however, we also need a geometric environment, which we 
introduce using presheaves.

A map among posets is said to be {\it open} 
iff it is open in the topological sense (posets can
be viewed as topological spaces whose open subsets are the downward
closed subsets); thus $f : Q\lora P$ is open iff it
is order-preserving and moreover satisfies the following condition
for all $q\in Q, p\in P$
$$
 p\leq f(q) ~\Rightarrow~ \exists q'\in Q~(q'\leq q ~\&~f(q')=p)~~. 
$$
Let us recall that open surjective maps are called p-morphisms
  in the standard non classical logics terminology. 

  Let $\Pzero$ be the category of finite rooted posets and open maps
between them; a \emph{presheaf} over $\Pzero$ is a contravariant
functor from $\Pzero$ to the category of sets and function, that is, a
functor $H:\Pzero^{op}\lora \bf Set$.  Let us recall what this means:
a functor $H:\Pzero^{op}\lora \bf Set$ associates to each finite
rooted poset $P$ a set $\FH{P}$; if $f : Q \lora P$ is an open map,
then we are also given a function $\FH{f} : \FH{P} \lora \FH{Q}$;
moreover, identities are sent to identities, while composition is
reversed, $H(g \circ f) = H(f) \circ H(g)$.

Our presheaves form a category whose objects are presheaves over
$\Pzero$ and whose maps are natural transformations; recall that a
natural transformation $\psi: H\lora H'$ is a collections of maps
$\psi_P:\FH{P}\lora \FH[H']{P}$ (indexed by the objects of $\Pzero$)
such that for every map $f:Q\lora P$ in $\Pzero$, we have
$\FH[H']{f}\circ \psi_P = \psi_Q\circ \FH{f}$. Throughout the paper,
we shall usually omit the subscript $P$ when referring to the
$P$-component $\psi_P$ of a natural transformation $\psi$.

The basic example of presheaf we need in the paper is described as follows.
 Let $L$ be a finite poset and let $h_L$ be the contravariant functor so defined:
 \begin{itemize}
  \item for a finite poset $P$, $h_L(P)$ is the set of all $L$-evaluations;
  \item for an open map $f:Q\lora P$, $h_L(f)$ takes $v: P\lora L$ to
    $v\circ f: Q\lora L$.
 \end{itemize}
The presheaf $h_L$ is actually a sheaf (for the canonical Grothendieck topology over $\Pzero$); we won't need this fact,\footnote{ 
The sheaf structure becomes essential for instance when one has to compute images - images are the categorical counterparts of second order quantifiers, 
see~\cite{GhilardiZawadowski2011}.
}
but 
we nevertheless call $h_L$ the \emph{sheaf of $L$-evaluations} (presheaves of the kind $h_L$, for some $L$, are called \emph{evaluation sheaves}).

Notice the following fact: if $\psi: h_L \lora h_{L'}$ is a natural
transformation, $v\in h_L(P)$ and $p\in P$, then
$\psi(v_p)= (\psi(v))_p$ (this is due to the fact that the inclusion
$\downarrow p \subseteq P$ is an open map, hence an arrow in
$\Pzero$); thus, we shall feel free to use the (non-ambiguous)
notation $\psi(v)_p$ to denote $\psi(v_p)= (\psi(v))_p$.

The notion of \emph{bounded bisimulation index} (\emph{\bindex}, for
short)\footnote{ This is called 'index' tout court
  in~\cite{GhilardiZawadowski2011}; here we used the word `index' for
  a different notion, since Section~\ref{sec:intro}.  }  takes
together structural and combinatorial aspects.  We say that a natural
transformation $\psi: h_L \lora h_{L'}$ \emph{has \bindex $n$} if, for
every $v:P\lora L$ and $ v': P'\lora L$, we have that $v\sim_n v'$
implies $\psi(v)\sim_0 \psi(v')$.

The following Proposition lists basic facts about \bindex{es}. In
particular, it ensures that natural transformations having a \bindex
compose.
\begin{proposition}
  \label{prop:indexDecreases}
  Let $\psi: h_L \lora h_{L'}$ have \bindex $n$; then it has also
  \bindex $m$ for every $m\geq n$. Moreover, for every $k\geq 0$, for
  every $v:P\lora L$ and $ v': P'\lora L$, we have that
  $v\sim_{n+k} v'$ implies $\psi(v)\sim_k \psi(v')$.
\end{proposition}

\begin{proof} Suppose that $\psi$ has \bindex $n$;
 we prove by induction on $k$ that
\[  \forall {v,v' } \; \; {\it if} \; \;  v
\sim_{n+k} v' \; \; {\it then} \; \; \psi(v) \sim_k \psi(v')
\hspace{20mm}  
{\rm (*)_k} \] 
For $k=0$, ${\rm (*)_k}$ is just the definition of $\psi$ having \bindex $n$.  
Suppose that ${\rm (*)_k}$ holds for some $k$. Let $v, v'$ be such that $v \sim_{n+k+1} v'$. 
We shall prove
that (let $P, P'$ be the domains of $v, v'$ respectively)
\[ \forall {p \in P} \; \exists {p' \in P'} \; \; \psi(v)_p \sim_k \psi(v')_{p'}  \]
(the converse statement is similar).
Fix $p \in P$.  Since $v \sim_{n+k+1} v'$, there is $p' \in P'$
such that $v_p \sim_{n+k} v'_{p'}$.  Using the inductive
assumption and the naturality of $\psi$, we obtain:
\[ \psi(v)_p = \psi(v_p) \sim_k \psi(v'_{p'}) = \psi(v')_{p'}  \]
as wanted.
\end{proof}

We are now ready to state duality theorems. As it is evident from the
discussion in Section~\ref{sec:classical}, it is sufficient to state a
duality for the category of finitely generated free Heyting algebras;
although it would not be difficult to give a duality for finitely
presented Heyting algebras, we just state a duality for the
intermediate category of Heyting algebras freely generated by a finite
bounded distributive lattice (this is quite simple to state and is
sufficient for proving \RT).

 \begin{theorem}\label{thm:duality}
   The category of Heyting algebras freely generated by a finite
   bounded distributive lattice is dual to the subcategory of
   presheaves over $\Pzero$ having as objects the evaluations sheaves
   and as arrows the natural transformations having a \bindex.
 \end{theorem}
 
 \noindent
 We present a full proof of the above Theorem in the next
   section.
 
 It is important to notice that in the subcategory mentioned in the
 above Theorem, products are computed as in the category of
 presheaves. This means that they are computed pointwise, like in the category of sets:
 in other words, we have that  $(h_L \times h_{L'})(P) = h_L(P) \times h_{L'}(P)$
 and $(h_L \times h_{L'})(f) = h_L(f) \times h_{L'}(f)$, for all $P$ and $f$.
  Notice moreover
   that $h_{L\times L'}(P)\simeq h_L(P) \times h_{L'}(P)$, so we have
   $h_{L\times L'}\simeq h_L \times h_{L'}$; in addition, 
  the two product
 projections  have \bindex 0.
   The situation strongly contrasts with other kind of dualities, 
   see \cite{Esa74} for example, 
   for which products are difficult to compute. 
   The ease by which products are computed
   might be seen as the principal reason for tackling a proof of
   \RT by means of sheaf duality.

  As a final information, we need to identify the dual of the free Heyting algebra on one generator:
  
  \begin{proposition}
    The dual of the free Heyting algebra on one generator is
    $h_{\two}$, where $\two$ is the two-element poset $\set{ 0, 1}$
    with $1\leq 0$.
 \end{proposition}
 Indeed, we shall see in the next Section that $h_{\two}$ is dual to
 the \Ha freely generated by the distributive lattice $\cD(\two)$, the
 lattice of downsets of the chain $\two$. Since $\cD(\two)$---which is
 a three element chain---is the free distributive (bounded) lattice on
 one generator, a standard argument proves that the \Ha freely
 generated by the distributive lattice $\cD(\two)$ is itself free on
 one generator.

\section{Proof of the Duality Theorem}
\label{sec:proofOfDuality}
\label{sec:dualityproof}
We present in this Section a proof of Theorem~\ref{thm:duality}. The
reader interested in \RT might wish to proceed directly to
Section~\ref{sec:finite}.  
  While the material in this Section is adapted
  from~\cite{GhilardiZawadowski97}, Theorem~\ref{thm:Vduality},
  generalizing the duality to some subvarieties of \Ha{s}, is new. 

\smallskip

With each $L$-evaluation $u:P\ra L$ and each $n\in\omega$ we associate
the set $Type_{n}(u)$\label{type} of $\sim_{n}$-equivalence classes,
$Type_{n}(u)\eqdef\set{ [u_p]_{n} \mid p\in P}$---where we
  recall that $[u_p]_{n}$ denotes the $\sim_{n}$-equivalence class of
$u_p$.  
An important, although simple, fact is given by the following
proposition: \vskip 2mm
\begin{proposition}\label{p41.2} 
  For a finite poset $L$ and $n\in \omega$,
  there are only finitely many equivalence classes of $L$-evaluations
  with respect to $\sim_n$.
\end{proposition}
\begin{proof}
  This is evident for $n=0$. For $n>0$, we argue by induction as
  follows. By Proposition \ref{p41.1}, we have that $u\sim_n v$ iff
  $Type_{n-1}(u)=Type_{n-1}(v)$, hence there cannot be more non
  $\sim_n$-equivalent $L$-evaluations than sets of $\sim_{n-1}$
  equivalence classes.
\end{proof}

Let ${\cal S}(h_L)$ be the set of subpresheaves $S$ of
$h_L$ satisfying the following condition for some $n\geq 0$
\begin{equation} \label{forall1}
\forall u:P \ra L,~\forall v:Q\ra L~ (u\in S_P~\&~u\sim_n
v~\Rightarrow~ v\in S_Q)\,.
\end{equation}
When the condition above holds, we say that $n$ is a \emph{\bindex}
for $S$.  
Notice that the choice of the naming \bindex is consistent with
  the one used in the previous Section. Indeed, for
  $S \subseteq h_{L}$, let $\chi : h_{L} \rto h_{\two}$ be defined by
  $\chi_{P}(u)(p) = 1$ if and only if $u_{p} \in S_{\da p}$. 
  If $n$ is a \bindex for $S$, then $\chi$ is a natural transformation
  and $n$ is a \bindex for $\chi$.  Indeed, $h_{\two}$ is a subobject
  classifier for subpresheaves that are sheaves for the canonical
  topology, see p.95 of \cite{GhilardiZawadowski2011}.

The definition of ${\cal S}(h_L)$ can be given 
in a slightly different
way by introducing the relations $\leq_n$.  We put:
\begin{enumerate}
 \item[{\rm (i)}] $v\leq_{0}u$ {\it iff}  $v(\rho)\leq u(\rho)$;
 \item[{\rm (ii)}] $v\leq_{n+1}u$ {\it iff}
   $\forall q\in Q\;\exists p\in P~(v_p\sim_n u_q)$.
\end{enumerate}

\begin{lemma}
  \label{lemma:equivforall1forall2} ${\cal S}(h_L)$ can be equivalently defined as the set of subpresheaves $S$ of
$h_L$ satisfying the following condition for some $n\geq 0$
  \begin{equation} \label{forall2} \forall u:P \ra L, ~\forall v:Q\ra
    L~ (u\in S_P~\&~v\leq_n u~\Rightarrow~ v\in S_Q)\,.
  \end{equation}  
\end{lemma}
\begin{proof}
  Let us call (for the time being) ${\cal S}'(h_L)$ the set of
  subpresheaves $S$ of $h_L$ satisfying
  condition~\eqref{forall2}. Clearly,
  ${\cal S}'(h_L)\subseteq {\cal S}(h_L)$.  For the converse, take
  $S\in {\cal S}(h_L)$ having \bindex $n$; in order to show that
  $S\in {\cal S}'(h_L)$, we show that it satisfies~\eqref{forall2} for
  $n+1$. Let in fact $u,v$ be such that $u\in S_P$ and
  $v\leq_{n+1} u$.  Then (considering the root of the domain of $v$)
  we know that there is $p\in P$ such that $v\sim_n u_p$; since $S$ is
  a subpresheaf of $h_{L}$, $u_p\in S_{\downarrow p}$ and finally
  $v\in S_Q$ because $n$ is a \bindex for $S$.
\end{proof}
Whenever a subpresheaf $S$ satisfies condition~\eqref{forall2}
relative to $n$, we say that $S$ has \blindex $n$.  Notice that,
  from these definitions, if $S$ has \blindex $n$, then it also has
  \bindex $n$, and if $S$ has \bindex $n$, then it has \blindex
  $n+1$.
It can be shown that $S$ has a \bindex $n$ iff it has \blindex $n$:
however, we won't use this result, since it depends on a construction
(the `grafting construction', see p.77 of~\cite{GhilardiZawadowski97})
which is not available if we  move from
the variety of \Ha{s} to one of its subvarieties.
Depending on the context, we shall make  use or not of the equivalent
definition for $\cS(h_{L})$ supplied by Lemma~\ref{lemma:equivforall1forall2}.

Let, for every
$u : P \rto L$ and $n \in \omega$,
\begin{align*}
  (\downarrow_{n}u)_Q & \eqdef\set{ v:Q\ra L \mid v\leq_{n} u }\,.
\end{align*}
The next Lemma is an immediate consequence of
Lemma~\ref{lemma:equivforall1forall2}.

\begin{lemma}
  $\da_{n} u$ is the least subpresheaf of $h_{L}$ having \blindex $n$
  such that $u \in F(P)$. A subpresheaf of $S$ of $h_{L}$ has \blindex
  $n$ if and only if, for each $u:P \ra L$ with $u \in S_P$,
  $\da_{n} u \subseteq S$.
\end{lemma}
In particular $\da_{n} u \in \cS(h_{L})$, for each $u : P \rto L$.  Notice
that the map
\begin{align*}
  [u]_{n} & \,\mapsto \,\da_{n} u\,,
\end{align*}
is well defined (actually, it is also injective) and so, by
Proposition~\ref{p41.2} and for fixed $n \in \omega$, there exists
only a finite number of presheaves of the form $\da_{n} u$. Since
\begin{align*}
  S & = \bigcup_{u \in S_{P}} \da_{n} u\,,
\end{align*}
when $S \in \cS(h_{L})$ has \blindex $n$, it follows that:
\begin{lemma}
  \label{lemma:finiteUnion}
  Every $S \in \cS(h_{L})$ of \blindex $n$ is a finite union of elements of
  the form $\da_{n} u$.
\end{lemma}

\medskip

Recall that $Sub(h_L)$ denotes the Heyting algebra of subpresheaves of
$h_L$.  
\begin{proposition}
  $\cS(h_L)$ is a sub-Heyting algebra of $Sub(h_L)$.
\end{proposition}
\begin{proof}
  It is easily seen that if $S$ and $T$ have \bindex $n$, then both
  $S \cap T$ and $S \cup T$ have \bindex $n$.  Next, consider the
  standard characterization of implication in subpresheaves:
  \begin{align*}
    (S\rightarrow T)_P & =\set{ u \in (h_L)_P \mid \forall h : Q \ra P
      ~(u \circ h\in {S}_{Q} ~\Ra~ u \circ h\in {T}_{Q})}\,.
  \end{align*}
  Notice that, for any $h:Q\ra P$ and $u \in (h_L)_P$, we have that $u\circ h\sim_{\infty} u_p$, where $p$ is $h(\Root{Q})$;
  as a consequence, since every $U\in \cS(h_L)$ has a \bindex, we have $u\circ h \in U_Q$ iff $u_p\in U_{\da p}$ for every $U\in \cS(h_L)$. 
  Thus, the following is an equivalent
  description of the implication
  \begin{align}
    \label{exists2}
    (S\rightarrow T)_P & =\set{ u\in (h_L)_P \mid \forall p\in P
      ~(u_p\in S_{\da p} ~\Ra~ u_p\in T_{\da p})}\,.
  \end{align}
  From this description it easily follows that if $S, T\in Sub(h_L)$
  have \bindex $n$, then $S\rightarrow T$ has \bindex $n+1$.
\end{proof}
\medskip

Let $\cD(L)$ denote the distributive lattice of downward closed
subsets of $L$ and recall that $\cD(L)$ is the Birkhoff dual of
the poset $L$, see \cite{birkhoff1937,DP}. Notice that there is a lattice embedding
$\iota_L:\cD (L)\ra \cS(h_L)$ associating with a downward closed
subset $d$ of $L$, the subpresheaf
\begin{align*}
  \iota_L(d)_P & \eqdef \set{ u: P\ra L \mid u(\rho(P))\in d}\,.
\end{align*}
Thus, for $p \in P$ and $u \in h_{L}(P)$, we have
$u_{p} \in \iota_{L}(d)_{\downarrow p}$ iff $u(p) \in d$.

We shall prove
that $\cS(h_L)$ is the free Heyting algebra generated by the
finite distributive lattice $\cD( L)$ with $\iota_L$ as the canonical
embedding.

\begin{lemma}\label{l42.2} 
The image of $\iota_L$ generates $\cS(h_L)$ as a Heyting
algebra.
\end{lemma}
\begin{proof}
  Clearly, the elements of $\cS(h_L)$ having \blindex 0 are exactly the
  elements of the image of $\iota_L$. Now consider an element having
  \blindex $n+1$; by Lemma~\ref{lemma:finiteUnion}, it is a finite
  union of elements of the kind $(\downarrow_{n+1}u)$. 
  We can express such elements in terms of elements having \blindex $n$
  as follows:
  \begin{equation}\label{nform}
    (\downarrow_{n+1} u)~=~ \bigcap_{p\in
        dom(u), \,v\not\sim_n u_p} (\,(\downarrow_n
    v)\rightarrow\bigcup_{v\not\leq_n w}
    (\downarrow_n w)\,)\,.
  \end{equation}
  Notice that all intersections and unions involved in the above
  formula are finite. Indeed, we have already observed that there are
  only finitely many elements of the kind $\downarrow_n w$. Moreover,
  if $v_1\sim_n v_2$, then $\da_{n} v_{1} = \da_{n} v_{2}$ and also
  $v_{1} \leq_{n} w$ if and only if
    $v_{2} \leq_{n} w$. As a consequence,
  $(\downarrow_n v_1)\rightarrow\bigcup_{v_1\not\leq_n w}
  (\downarrow_n w)$ equals
  $(\downarrow_n v_2)\rightarrow\bigcup_{v_2\not\leq_n w}
  (\downarrow_n w)$.
  \smallskip

Let us verify equation~\eqref{nform}. Suppose that $z\leq_{n+1}u$ and
let $v : Q \ra L$ be arbitrary with the
property 
that $v\not\sim_n u_p$, for every point $p$ in the domain of $u$. We
show that
$z\in (\downarrow_n v)\rightarrow \bigcup_{\lbrace w\mid v\not\leq_n
  w\rbrace} (\downarrow_n w)$ using~\eqref{exists2}.
Let $q$ be a point in the domain of $z$ such that $z_q\leq_n v$. From
$z\leq_{n+1}u$ we conclude that there exists $p$ such that
$z_q\sim_n u_p$. Consequently, $v\not\leq_n z_q$, otherwise
$z_q\sim_n v$ and so $v\sim_n u_p$, contradicting the choice of
$v$. Therefore
$z_q\in \, \da_n z_{q} \subseteq \bigcup_{v\not\leq_n w} (\da_n
w)$. Vice versa, suppose that $z\not\leq_{n+1}u$. It follows that
there is a point $q$ in the domain of $z$ such, that for every point
$p$ in the domain of $u$, $z_q\not\sim_n u_p$. We check that
$z\not\in (\da_n z_q)\ra \bigcup_{z_q\not\leq_n w
}\da_n(w)$. This is clear as $z_q \in (\da_n z_q)$ and
$z_q\not\in\bigcup_{z_q\not\leq_n w} (\da_nw)$.

This proves equation~\eqref{nform} and ends the proof of the Lemma.
\end{proof}

\medskip

The following statement is an immediate consequence of the finite
model property:
\begin{lemma}
  \label{lemma:productOfFinite}
  Every finitely presented \Ha embeds into  a product of finite \Ha{s}.  
\end{lemma}

Recall that, by Birkhoff duality, monotone maps $f:M \ra L$ between
finite posets $M,L$ bijectively (and naturally) correspond to
bound-preserving lattice homomorphism
$f^{-1} = \cD(f) : \cD(L) \rto \cD(M)$.  Therefore, for $f:M \ra L$,
we define a map
\begin{align}
  \notag
  ev_f &:  \cS(h_L) \longrightarrow \cD(M)
  \intertext{by putting, for
    $X \in \cS(h_L)$,}
  \label{eq:defEv}
  ev_f(X) & \eqdef \set{ p \in M \mid f_p \in X_{\da p}} \,.
\end{align}
\begin{proposition}
  \label{prop:freeness}
  The map $ev_{f}$ is a \Ha morphism and makes the following diagram
  commute:
  \begin{center}
    \settriparms[1`1`1;550] \qtriangle[\cD(L)`\cS(h_L)`{\cal
      D}(M);\iota_L`\cD(f)`ev_{f}]
  \end{center}
  Consequently, $\cS(h_L)$, together with $\iota_L$ as the canonical
  embedding, is a free \Ha generated by the finite distributive
  lattice $\cD(L)$.
\end{proposition}
\begin{proof}
  Let us verify first that the above diagram commutes. For each $d
  \in \cD(L)$,
  \begin{align*}
    ev_{f}(i_{L}(d))
    & = \set{p \in M  \mid f_{p} \in i_{L}(d)_{\da p}} \\
    & = \set{p \in M \mid f_{p}(p) \in d} \\
    & = \set{p \in M \mid f(p) \in d} = \cD(f)(d)\,.
\end{align*}
To see that $ev_{f}$ is a \Ha homomorphism, we have, for example,
\begin{align*}
  ev_{f}(S \rightarrow T)
  & = \set{p \in M  \mid f_{p} \in (S \rightarrow T)_{\downarrow p}} \\
  & = \set{p \in M \mid \forall q \in \downarrow p \,(f_{pq} \in
    S_{\downarrow q}
    \Rightarrow f_{pq} \in T_{\downarrow q})}  \\
  & = \set{p \in M \mid \forall q \leq p \,(f_{q} \in
    S_{\downarrow q} \Rightarrow f_{q} \in T_{\downarrow q})\,} \\
  &=
  ev_{f}(S) \rightarrow ev_{f}(T)\,.
\end{align*}
Notice also that, in view of Lemma~\ref{l42.2}, $ev_{f}$ is the unique
\Ha morphism $g : \cS(h_{L}) \rto \cD(M)$ with the property that
$g\circ i_{L} = \cD(f)$. Therefore, we have argued that every bounded
lattice morphism $g = \cD(f) : \cD(L) \ra \cD(M)$, where $M$ is a
finite poset, extends uniquely to the \Ha morphism
$ev_{f} : \cS(h_{L}) \ra \cD(M)$.
By a standard argument, the same universal property holds with respect
to the bound-preserving lattice homomorphisms $g : \cD(L) \ra H$,
where now $H$ is a sub-\Ha of a product of finite \Ha{s} of the form
$\cD(M)$. In particular, using Lemma~\ref{lemma:productOfFinite}, we
can take $(H,g)$ to be $(F,\eta)$, the free \Ha algebra generated by
the distributive lattice $\cD(M)$. Then, by combining the universal
properties of $(\iota_{L},\cS(h_{L}))$ and of $(F,\eta)$, it follows
that $\cS(h_{L})$ and $F$ are isomorphic.
\end{proof}

\medskip

Let $\HD$ be the category of Heyting algebras freely generated by
a finite distributive lattice and let $\MH$ be the subcategory
of presheaves over $\Pzero$ having as objects the evaluations sheaves
and as arrows the natural transformations having a \bindex.
We want to show that
$\HD$ is dual to $\MH$.

We define the following functor $\TH$: \vskip 4mm
\begin{center}
  \begin{tikzcd}
    \MH \arrow[rrr,"\TH"]&&& \;\;\HD^{op} \\
    h_{L} \arrow[dd,"f"']&&& \cS(h_{L}) \\
    {}\arrow[rrr, mapsto]&&{}&{}\\
    h_{M} &&& \cS(h_{M}) \arrow[uu,"f^{-1} = \TH(f)"']
  \end{tikzcd}
\end{center}
where $\cS(h_{N})$ is as in Lemma \ref{l42.2} and $f^{-1}$ is the
inverse image function.  \vskip 5mm\noindent
\begin{lemma}\label{l43.2} 
\begin{description}
\item[{\rm (i)}]  ${\bf T_H}$  is a well defined functor.
\item[{\rm (ii)}]  ${\bf T_H}$  is essentially surjective.
\end{description}
\end{lemma}
\begin{proof}
  By Proposition~\ref{prop:freeness},
$\cS(h_L)$ is a
Heyting algebra freely generated by a finite distributive lattice and
every such Heyting algebra is isomorphic to one of that form. Hence
${\bf T_H}$ is well defined on objects and essentially
surjective. Clearly ${\bf T_H}$ preserves compositions and
identities. We need to show that for any subpresheaf $D$ of $h_M$ with
a b-index, $f^{-1}(D)$ has a b-index and that $f^{-1}$ is a Heyting
algebra morphism. The latter follows from the fact that
$f^{-1}: Sub(h_M) \longrightarrow Sub(h_L)$ is a Heyting algebra
morphism and that $\cS(h_L)$, $\cS(h_M)$ are sub-Heyting
algebras of the Heyting algebras of subpresheaves
$Sub(h_L), Sub(h_M)$, respectively.

Let $n$ be a \bindex of $D$ and $m$ a \bindex of $f$.  We shall show
that $f^{-1}(D)$ has \bindex $n+m$.  Let $v~\in~f^{-1}(D)$ and $v'$ be
such that $v~\sim_{n+m}~v'$. By Proposition~\ref{prop:indexDecreases},
$f(v) \sim_{n} f(v')$.  Since $f(v) \in D$ and $D$ is $\sim_n$-closed,
it follows that $f(v') \in D$ and then $v' \in f^{-1}(D)$.
\end{proof}

\medskip

Recall from Proposition~\ref{prop:freeness} that, for $u :P \ra L$,
$ev_{u}$ is the unique \Ha morphism $\cS(h_L) \longrightarrow \cD(P)$
such that $\cD(u) = ev_{u} \circ \iota_{L}$.
Conversely, given a \Ha morphism
$\alpha : \cS(h_L) \longrightarrow \cD(P)$, we define an
$L$-evaluation
\[ \overline{\alpha} :P \ra L \] as the dual of the distributive
lattice morphism $\alpha \circ \iota_{L}$. By the definition of
$\overline{\alpha}$, the diagram
\begin{center}
  \begin{tikzcd}
    \cD(L) \arrow[rr,"\iota_{L}"]
    \arrow[rrd,"\cD(\overline{\alpha})"']
    && \cS(h_{L}) \arrow[d,"\alpha"] \\
    && \cD(L)
  \end{tikzcd}
\end{center}
commutes and therefore, by the universal property of $ev$, we
deduce the following relation:
\begin{align}
  \label{eq:evBarAlpha}
  ev_{\overline{\alpha}} & = \alpha\,.
\end{align}

The two maps
\[ \alpha \mapsto \overline{\alpha} \;\;\;\;\;\; u \mapsto ev_u\,, \]
yield a bijective correspondence between the \Ha morphisms
$\alpha : \cS(h_L) \ra \cD(P)$ and the $L$-evaluations
$u \in h_{L}(P)$ which is natural in $P$.  This immediately follows
from the chain of natural isomorphisms
\begin{align*}
  \Cat{POS}(P,L) & \simeq \Cat{DLATT}(\cD(L),\cD(P)) \simeq
  \Cat{HA}(\cS(h_{L}),\cD(P))\,,
\end{align*}
where the first natural isomorphism is Birkhoff duality between the
category of finite posets and the category of finite distributive
lattices, and the second is by freeness of $\cS(h_{L})$,
Proposition~\ref{prop:freeness}.
\color{black}

Let $h_L, h_M$ be objects of $\HD$, $\mu :\cS (h_L) \ra \cS (h_M)$ be
a morphism of Heyting algebras.  For each $P\in \Pzero$, we define 
\begin{align*}
  \mu^*_P & :\;h_M(P) \longrightarrow h_L(P)
  \intertext{as follows:}
  \mu^*_P(u) & \eqdef \overline{ev_u\circ\mu}\,, \quad \text{for each $u\in h_M(P)$.}
\end{align*}
Note that by the above correspondence, it is immediate that
$\mu^*_P(u)\in h_L(P)$ and that $\mu^*:h_M \ra h_L$ is a natural
transformation. Moreover

\vskip 2mm \noindent
\begin{proposition}\label{p43.3new}
With the notation as above, we have
\begin{enumerate}[{\rm (i)}]
\item 
  $\mu_{P}^*(u)\in X_{P}$ iff $u\in \mu(X)_{P}$, for $X\in\cS (h_L)$ and $u\in h_M(P)$;
   \label{it:iff}
 \item $\mu^*: h_M \ra h_L$ is a morphism in $\bf M_H$;
   \label{it:morphism}
 \item $\mu = (\mu^*)^{-1} = \TH(\mu^{*})$;
   \label{it:full}
 \item $f = (f^{-1})^*$, for any morphism $f: h_M \longrightarrow h_L$ in $\bf M_H$.
  \label{it:faithful}
\end{enumerate}
\end{proposition}
\begin{proof}
  Ad \eqref{it:iff}. Observe that, for any $v \in h_{M}(Q)$ and
  $Y \in \cS(h_L)$, $ev_{v}(Y) = Q$ iff, for all $q \in Q$,
  $v_{q} \in Y_{\da q}$, iff $v_{\rho(Q)} \in Y_{\da \rho(Q)}$ iff
  $v \in Y_{Q}$.
  Recall now---see equation~\eqref{eq:evBarAlpha}---that, for any
  $\alpha: \cS(h_{L}) \ra \cD(M)$, $ev_{\overline{\alpha}} = \alpha$
  and so, in particular,
  $ev_{\mu_{P}^*(u)} = ev_{\,\overline{ev_u\circ \mu}} = ev_u\circ
  \mu$.  Therefore, for $X\in \cS (h_L)$ and $u:P\ra M$, we have
  $ev_{\mu_{P}^*(u)}(X)=ev_u(\mu(X))$ and therefore, according to the
  previous observation, we have $\mu_{P}^*(u)\in X$ iff
  $ev_{\mu_{P}^*(u)}(X)=P$ iff $ev_u(\mu(X))=P$ iff
  $u\in \mu (X)_{P}$.

  Ad~\eqref{it:morphism}. We need to show that the transformation $\mu^*$ has a
  \bindex. Let $n\in\omega$ be the maximum of the \bindex{es} 
  of sets of $\mu (X)$ where $X$ is of the kind $\iota_L(d)$ for some
  $d\in \cD(L)$.  Notice that there are only finitely many such $X$'s
  and, moreover, for $w, w'\in h_L$ we have $w\sim_0 w'$ iff $w, w'$
  belong to the same such $X$'s.
  For any
  $u \in h_{M}(P)$ and $v \in h_{M}(Q)$ such that $u\sim_n v$ and for
  $X= \iota_L(d)$, we have
  $$
  \begin{equivJustif}
    \mu^*_{Q}(v) \in X_{Q}
    \iffJustif[by~\eqref{it:iff}] 
    v \in \mu(X)_{Q} 
    \iffJustif
    u \in \mu(X)_{P} 
    \iffJustif
    \mu^*_{P}(u) \in X_{P} 
  \end{equivJustif}
  $$
  where the horizontal lines above stand for logical
    equivalences. 
  Thus $\mu_{Q}^*(v) \sim_0 \mu_{P}^*(u)$ and $\mu^*$ has \bindex $n$.
  
  Ad~\eqref{it:full}. Using~\eqref{it:iff}, we have,  for any $X\in\cS(h_L)$ and $v\in
  h_M(P)$,
  $$
  \begin{equivJustif} 
    v\in (\mu_{P}^*)^{-1}(X) 
    \iffJustif 
    \mu_{P}^*(v) \in X 
    \iffJustif 
    v \in \mu(X)_{P}
  \end{equivJustif}
  $$
  i.e. $\mu=(\mu^*)^{-1}$.

  Ad~\eqref{it:faithful}. Let $v \in h_M(P)$, $p\in P$ and $d\in {\mathcal D}
  (L)$. Then, we have
  $$
  \begin{equivJustif}
    (f^{-1})_{P}^*(v)(p)\in d 
    \iffJustif[by the definition of $\iota_{L}$,] 
    (f^{-1})_{\da p}^*(v_p) \in \iota_L(d)_{\da p} 
    \iffJustif[using~\eqref{it:iff},] 
    v_p\in (f^{-1}(\iota_L(d)))_{\da p} 
    \iffJustif 
    f_{\da p}(v_p) 
    \in \iota_L(d)_{\da p} 
    \iffJustif 
    f_P(v)(p)\in d 
  \end{equivJustif}
  $$
  Since $P$, $v$, $p$ and $d$ were arbitrary  $f = (f^{-1})^*$.
\end{proof}

\medskip

Thus we have :
\begin{theorem}[Duality Theorem]
  \label{t43.5} 
  The functor $\TH : \MH \longrightarrow \HD^{op}$ is an
  equivalence of categories.
\end{theorem}
\begin{proof}
  Lemma \ref{l43.2} shows that $\TH$ is a functor which is
  essentially surjective and by Lemma
  \ref{p43.3new}(\ref{it:morphism}--\ref{it:faithful}) $\TH$ is
  full and faithful, i.e. $\TH$ is an equivalence of categories.
\end{proof}

For some applications in Section~\ref{sec:bounds}, we shall need a
duality theorem for some subvarieties. Call a variety $\bf V$ of
Heyting algebras \emph{finitely approximable} if every finitely
generated free $V$-algebra embeds into a product of finite
$\bf V$-algebras.

We can extend the above duality Theorem to finitely approximable subvarieties as follows. Take one such subvariety $\bf V$ and let $\Pzero^{\bf V}$
be the category of finite rooted posets $P$ such that $\cD(P)\in \bf V$.
Let $\HD^{\bf V}$ be the category of $\bf V$-algebras freely generated by
a finite distributive lattice and let $\MH^{\!\bf V}$ be the subcategory
of presheaves over $\Pzero^{\bf V}$ having as objects the evaluations sheaves
and as arrows the natural transformations having a \bindex. We have:

\begin{theorem}[Duality Theorem for Finitely Approximable Subvarieties]
  \label{thm:Vduality} For every finitely approximable variety $\bf V$ of Heyting algebras,
  $\HD^{\bf V}$ is dual to $\MH^{\!\bf V}$. 
\end{theorem}
\begin{proof} 
  By reading back
  the proof of Theorem~\ref{t43.5}, it is immediately realized that
  Lemma~\ref{lemma:productOfFinite} is the only specific fact on
  Heyting algebras we used.  When this Lemma is replaced by the
  assumption that $\bf V$ is finitely approximable, the same chain of
  arguments yields a proof of Theorem~\ref{thm:Vduality}.
\end{proof}

\section{Indexes and Periods over Finite Models}
\label{sec:finite}

Taking into consideration the algebraic reformulation from
Section~\ref{sec:classical} and the information from Section~\ref{sec:dualityproof}, we can prove
\RT for $(IPC)$ by showing that \emph{all natural transformations from
  $h_{L} \times h_{\bf 2}$ into itself, commuting over the first
  projection $\pi_0$ and having a b-index, are ultimately periodic
  with period 2}.  Spelling this out, this means the following. Fix a
finite poset $L$ and a natural transformation
$\psi: h_{L} \times h_{\bf 2}\lora h_{L} \times h_{\bf 2}$ having a
b-index such that the diagram
  \begin{center}
\resetparms
\settriparms[1`1`1;400] \Vtriangle[ h_{L}\times h_{\bf 2}` h_{L}\times h_{\bf 2} 
`h_{L} ;\psi`\pi_0`\pi_0 ]
\end{center}
\noindent
commutes; we have to find an $N$ such that $\psi^{N+2}=\psi^N$,
according to the dual reformulation of~\eqref{eq:Ruitenburgtheorem}.

From the commutativity of the above triangle, we can decompose $\psi$
as $\psi = \langle \pi_{0},\chi\rangle$, were both
$\pi_0 : h_{L} \times h_{\bf 2} \lora h_{L}$ and
$\chi : h_{L} \times h_{\bf 2} \lora h_{\bf
    2}$
  have a b-index; we assume that $n\geq 1$ is a b-index for both of
  them. \emph{We let such $\psi = \langle \pi_{0},\chi\rangle$ and $n$
    be fixed for the rest of the paper.}

Notice that for $\vu \in h_{L}(P)\times  h_{\bf 2}(P)$, we have
$$
{\psi^{k}}\vu= (v,u_k)
$$
where we put 
\begin{equation}\label{eq:def_un}
 u_0\eqdef u~~{\rm  and}~~ u_{{k}+1}\eqdef\chi(v,
 u_{{k}})
 \,.
\end{equation}
Since $P$  and $L$ are finite, it is clear that the sequence $\set{\psi^{k}\vu \mid k \geq 0}$ 
(and obviously also the sequence $\set{u_{k}\mid k \geq 0}$) must become ultimately periodic.

We show in this section that, for each finite set $P$ and for each
$\vu \in h_{L}(P)$, the period of the sequence
$\set{\psi^{k}\vu \mid k \geq 0}$ has $2$ as an upper bound, whereas
the index of $\set{\psi^{k}\vu \mid k \geq 0}$ can be bounded by the
maximum length of the chains in the finite poset $P$ (in the next
section, we shall bound such an index independently on $P$, thus
proving \RT).

Call $\vu \in h_{L}(P)$ \emph{2-periodic} (or just
\emph{periodic}\footnote{From now on, `periodic' will mean
  `2-periodic', i.e. `periodic with period 2'.}) iff we have
$\psi^2\vu=\vu$; a point $q\in P$ is similarly said periodic in
  $\vu$ iff $\vu_q$ is periodic. We shall only say that $p$ is
  periodic if an evaluation is given and understood from the
  context. We call a point \emph{\nonperiodic} if it is not periodic
  (w.r.t. a given evaluation).
  
\begin{lemma}\label{lem:period}
  Let $\vu \in h_{L}(P)$ and $p\in P$ be such that all $q\in P$,
  $q< p$, are periodic. Then either $\vu_p$ is periodic or $\psi\vu_p$
  is periodic.  Moreover, if $\vu_p$ is \nonperiodic and
  $u_0(p)=u(p)=1$, then $u_1(p)=\chi(u,v)(p)=0$.
\end{lemma}

\begin{proof} We work by induction on the height of $p$ (i.e. on the maximum $\leq$-chain starting with $p$ in $P$).
  If the height of $p$ is $1$, then the argument is the same as in the
  classical logic case (see Section~\ref{sec:classical}).

  If the height is greater than one, then we need a simple
  combinatorial check about the possible cases that might arise.
  Recalling the above definition~\eqref{eq:def_un} of the
  $\bf 2$-evaluations
  $u_n$, 
  the induction hypothesis tells us that there
  is $M$  big enough so that so for all $k \geq M$ and
  $q < p$, $(u_{k+2})_{q} = (u_{k})_{q}$.

Let $\ddownset p = \set{q \in P \mid q < p}$.  We shall represent
$(u_{k})_{p}$ as a pair $\typeOfP{a_{k}}{x_{k}}$, where
$a_{k} = u_{k}(p)$ and $x_{k}$ is the restriction of
$(u_{k})_{p}$ to $\ddownset p$.
  
Let us start by considering a first repeat $(i,j)$ of the sequence
$\set{a_{M +k}}_{k\geq 0}$ - that is $i$ is the smallest $i$ such that
there is $j > 0$ such that $a_{M +i+j}= a_{M +i}$
and $j$ is the smallest such $j$. Since the $a_{M +n}$ can only take
value 0 or 1, we must have $i+j \leq 2$.  We show that the sequence
$\set{(u_{M +k})_p}_{k\geq 0}$ has first repeat taken from
  $$
  (0,1), (0,2), (1,1), (1,2)\,.ùì
  $$ 
This shall imply in the first two cases
  that $\vu_p$ is periodic or, in the last two cases, that $\psi\vu_p$
  is periodic.  To our goal, let $x = x_{M}$ and $y = x_{M +1}$
(recall that we do now know whether $x = y$).

  Notice that, if $j = 2$, then $i = 0$ and a first repeat for
  $\set{(u_{k})_p}_{k \geq M}$, is $(0,2)$, as in the diagram below
  \begin{align*}
    \typeOfP{a}{x}\typeOfP{b}{y}\typeOfP{a}{x}\,.
  \end{align*}

  Therefore, let us assume $j = 1$ (so $i \in \set{0,1}$). Consider
  firstly $i = 0$:
  \begin{align*}
    \typeOfP{a}{x}\typeOfP{a}{y}\typeOfP{c}{x}\typeOfP{d}{y}
  \end{align*}
  If $x = y$, then we have a repeat at $(0,1)$. Also, if $a =1$, then
  the mappings $x$ and $y$ are uniformly $1$ (since evaluations
    are order-preserving maps and we have $1\leq 0$ in
    $\bf 2$):
  again, $x = y$ and
  $(0,1)$ is a repeat.

  So let us assume $x \neq y$ and $a = 0$. If $c = a$, then we have
  the repeat $(0,2)$ as above.  Otherwise $c = 1$, so $x = 1$.  We
  cannot have $d = 1$, otherwise $1 = x = y$. Thus
  $d = 0 = a$, and the repeat is $(1,2)$.

  Finally, consider $i = 1$ (so $a \neq b$ and $j = 1$):
  \begin{align*}
    \typeOfP{a}{x}\typeOfP{b}{y}\typeOfP{b}{x}\typeOfP{d}{y}
  \end{align*}
  We have two  subcases: $b=1$ and $b=0$. If $b=1$, then $a=0$ and $x = 1 = y$: we have a repeat at $(1,1)$.
  
  In the last subcase, we have $b=0$, $a=1$ and now if $d=0$ we have a
  repeat at $(1,2)$ and if $d=1$ we have a repeat $(1,1)$ (because
  $d = a = 1$ implies $y = 1$ and $x =1$).

  The last statement of the Lemma is also obvious in view
    of the fact that if $a = b = 1$, then $x = y =1$, so $p$ is
    periodic.
\end{proof}

\begin{corollary}
  \label{cor:indexFromHeight}
  Let $N_P$ be the height of $P$; then $\psi^{N_P}\vu$ is periodic for
  all $\vu \in h_{L}(P)$.
\end{corollary}

\begin{proof}
 An easy induction on $N_P$, based on the previous Lemma.
\end{proof}

\section{Ranks}\label{sec:ranks}

Ranks (already introduced in~\cite{fine}) are a powerful tool
that goes hand in hand with bounded
bisimulations; in our context the useful notion of rank is given
below.  Recall that $\psi = \langle \pi_{0},\chi\rangle$ and that
$n\geq 1$ is a b-index for $\psi$ and $\chi$.

Let $\vu \in h_{L}(P)$ be given.  The \emph{type} of a periodic point
$p\in P$ is  the pair of equivalence classes
\begin{equation}\label{eq:pairs}
 \langle [(v_p,u_p)]_{n-1}, [\psi(v_p,u_p)]_{n-1}\rangle.
\end{equation}
The \emph{rank} of a point $p$ (that we shall denote by
  $rk(p)$) is the cardinality of the set of distinct types of the
periodic points $q\leq p$.  Since $\sim_{n-1}$ is an equivalence
relation with finitely many equivalence classes, the rank cannot
exceed a positive number $R(L,n)$ (that can be computed in function of
$L, n$).

Clearly we have $rk(p)\geq rk(q)$ in case $p\geq q$.  Notice that an
application of $\psi$ does not decrease the rank of a point: this is
because the pairs~\eqref{eq:pairs} coming from a periodic point just
get swapped after applying $\psi$.  A non-periodic point $p\in P$ has
\emph{minimal rank} iff we have $rk(p)=rk(q)$ for all \nonperiodic
$q\leq p$.

\begin{lemma}\label{lem:minrank}
  Let $p\in P$ be a \nonperiodic point of minimal rank in
  $\vu \in h_{L}(P)$; suppose also that $\vu$ is constant
    on the set of all \nonperiodic points in
    $\downarrow p$.  Then we have
  $\psi^m\vu_{q_0}\sim_n\psi^m\vu_{q_1}$ for all $m\geq 0$ and for all
  \nonperiodic points $q_0, q_1\leq p$.
\end{lemma}
\begin{proof} We let $\Pi$ be the set of periodic points of $\vu$ that
  are in $\downarrow p$ and let $\Pi^c$ be
  $(\downarrow p)\setminus \Pi$ .  Let us first observe that for every
  $r\in \Pi^c$, we have
\begin{eqnarray*}
\smalleql{\{\langle 
 [(v_s,u_s)]_{n-1}, [\psi(v_s,u_s)]_{n-1}\rangle \mid s\leq r, ~s~\hbox{is  periodic}\}}
& \\
& =~~  \{\langle [(v_s,u_s)]_{n-1}, [\psi(v_s,u_s)]_{n-1}\rangle \mid s\leq p, ~s~\hbox{is periodic}\}
\end{eqnarray*}
(indeed the inclusion $\subseteq$ is because $r\leq p$ and the
inclusion $\supseteq$ is by the minimality of the rank of $p$).
Saying this in words, we have that ``for every periodic $s\leq p$
there is a periodic $s'\leq r$ such that
$(v_s,u_s)\sim_{n-1} (v_{s'},u_{s'})$ and
$\psi(v_s,u_s)\sim_{n-1} \psi(v_{s'},u_{s'})$''; also (by the
definition of 2-periodicity), ``for all $m\geq 0$, for every periodic
$s\leq p$ there is a periodic $s'\leq r$ such that
$\psi^m(v_s,u_s)\sim_{n-1} \psi^m(v_{s'},u_{s'})$''.
By letting both $q_0, q_1$ playing the role of
$r$, we get:
\begin{fact*}
  For every $m\geq 0$, for every
  $q_0, q_1\in \Pi^c$, for every periodic $s\leq q_0$ there is a
  periodic $s'\leq q_0$ such that
  $\psi^m(v_s,u_s)\sim_{n-1}
    \psi^m(v_{s'},u_{s'})$ (and vice versa).
\end{fact*}

We now prove the statement of the theorem by induction on $m$; take
two points $q_0, q_1\in \Pi^c$.

For $m=0$, $\vu_{q_0}\sim_n\vu_{q_1}$ is established as follows: as
long as Player 1  plays in $\Pi^c$, we know $\vu$ is constant so that
Player 2 can answer with an identical move still staying within $\Pi^c$;
as soon as it plays in $\Pi$, Player 2 uses the above Fact to win the game.

The inductive case
  $\psi^{m+1}\vu_{q_0}\sim_n\psi^{m+1}\vu_{q_1}$ is proved in the same
  way, using the Fact (which holds for the integer $m +1$) and
  observing that $\psi^{m+1}$ is constant on $\Pi^{c}$. The latter
  statement can be verified as follows: by the induction hypothesis we
  have $\psi^m\vu_{q}\sim_n\psi^m\vu_{q'}$, so we derive from
  Proposition~\ref{prop:indexDecreases}
  $\psi^{m+1}\vu_{q}\sim_0\psi^{m+1}\vu_{q'}$, for all
  $q, q'\in \Pi^c$; that is, $\psi^{m+1}$ is constant on $\Pi^{c}$. 
\end{proof}

\section{\RT}\label{sec:main}

We can finally prove: 
\begin{theorem} [\RT for IPC]\label{thm:main}
There is $N\geq 1$ such that we have $\psi^{N+2}=\psi^N$.
\end{theorem}

\begin{proof} Let $L$ be a finite poset and let $R\eqdef R(L,n)$ be
  the maximum rank for $n,L$ (see the previous section).  Below, for
  $e\in L$, we let $\vert e\vert$ be the height of $e$ in $L$,
  i.e. the maximum size of chains in $L$ whose maximum element is
  $e$; we let also $\vert L \vert$ be the maximum
  size of a chain in $L$.  We make an induction on
  natural numbers $l\geq 1$ and show the following: \emph{(for each
  $l \geq 1$) there is $N(l)$ such that for every $(v,u)$ and
  $p\in dom(v,u)$ such that $l\geq \vert v(p)\vert$, we have that
  $\psi^{N(l)}(v_p, u_p)$ is periodic.}     (It will turn out that $N(l)$ is
      $2R(l-1)+1$).
 Once this is proved,
  the statement
    of the Theorem shall be proved with
    $N= N(\vert L\vert)$.

  If $l=1$, it is easily seen that we can put $N(l)=1$ (this case is
  essentially the classical logic case).

  Pick a $p$ with $\vert v(p)\vert =l>1$; let $N_0$ be the maximum of
  the values $N(l_0)$ for $l_0< l$:\footnote{It is easily seen that we
    indeed have $N_{0} = N(l-1)$. }
    we show that we can take $N(l)$ to be $N_0+2R$.

    Firstly, let $(v, u_0) := \psi^{N_0}\vu$ so all $q$ with
    $\vert v(q)\vert < l$ are periodic in $(v,u_0)$.  After such iterations,
    suppose that $p$ is not yet periodic in $(v, u_0)$.  We let $r$ be
    the minimum rank of points $q\leq p$ which are not
    periodic 
    (all such points $q$ must be such that $v(q)=v(p)$); we show that
    after \emph{two iterations} of $\chi$, all points $p_0\leq p$
    having rank $r$ become periodic or increase their rank, thus
    causing the overall minimum rank below $p$ to increase: this means
    that after at most $2(R-r)\leq 2R$ iterations of $\psi$, all
    points below $p$ ($p$ itself included!)  become periodic
  (otherwise said, we take $R-r$ as the secondary parameter of our double induction).

Pick $p_0\leq p$ having minimal rank $r$;
 thus we have that all $q\leq p_0$
in $(v, u_0)$ are now either periodic or have the same rank and the same $v$-value  as $p_0$
(by the choice of $N_0$ above). Let us divide the  points of $\downarrow p_0$ into four subsets:
\begin{align*}
\Eper~ := &~ \set{ q \mid ~q ~{\rm is~periodic}}\\
E_0~~~\, := &~ \set{ q \mid~q\not\in \Eper~\&~\forall q'\leq q ~(q'\not\in \Eper \Rightarrow u_0(q')=0)} \\
E_1~~~\, := &~ \set{ q \mid~q\not\in \Eper~\&~\forall q'\leq q ~(q'\not\in \Eper \Rightarrow u_0(q')=1)}\\
E_{01}~~\, := &~ \set{ q \mid~ q'\not\in \Eper\cup E_1\cup E_0}\,.
\end{align*}
Let us define \emph{frontier point} a \nonperiodic point $f\leq p$ such that
all $q<f$ are periodic (clearly, a frontier point belongs to
$E_0\cup E_1$); by Lemma~\ref{lem:period}, all frontier points become
periodic after applying $\psi$.  Take a point
$q\in E_i$ and a frontier point $f$ below it;
 since $q$ also has
minimal rank and the hypotheses of Lemma~\ref{lem:minrank} are
satisfied for $\vu_q$, we have in particular that
$\psi^m(v,u_0)_{q'}=\psi^m(v,u_0)_f$
for all $m\geq 0$ and all \nonperiodic $q'\leq q$, and hence
$\psi(v,u_0)_q$ is periodic too.

Thus, if we apply $\psi$, we have that in $(v,u_1):=\psi(v,u_0)$ all
points in $\Eper\cup E_0\cup E_1$ become periodic, together with
possibly some points in $E_{01}$.  The latter points get in any case
  $u_1$-value equal to $0$. This can be seen as follows. If any
  such point gets $u_1$-value equal to $1$, then all points below it
  get the same $u_{1}$-value. Yet, by definition, these points are
  above some frontier point in $E_1$ and frontier points in $E_1$ get
  $u_1$-value $0$ by the second statement of Lemma~\ref{lem:period}.

If $p_0\in E_0$ has become
periodic, we are done; we are also done if the rank of $p_0$
increases, because this is precisely what we want.
If $p_0$ has not become periodic and its rank has not increased, then
now all the \nonperiodic points below $p_0$ in $(v,u_1)$ have
$u_1$-value $0$ (by the previous
  remark) and have the same rank as $p_0$. Thus, they are the set
$E_0$ computed in $(v, u_1)$ (instead of in $(v,u_0)$) and we know by
the same considerations as above that it is sufficient to apply $\psi$
once more to make them periodic.
\end{proof}
\begin{figure}
  \centering
\begin{tikzpicture}[xscale=0.7,yscale=.7]
\draw[very thick] (0,0) to [out=90,in=195] (1.7,4.5); 
\draw[very thick] (1.7,4.5) to [out=10,in=90] (3.4,0); 
\draw (0,0.6) --(3.4,1.8); 
\draw (.2,2.7) --(3.3,2.9); 
\draw (1.7,1.2) --(1.7,2.8); 
\node at (1.9,4.8) {$p_0$};
\node at (1.9,4.5) {$\bullet$};
\node at (1.0,1.9) {$E_0$};
\node at (2.5,1.9) {$E_1$};
\node at (1.7,3.4) {$E_{01}$};
\node at (1.7,0.4) {$E_{per}$};

\node at (4.9,2.7) {$\psi$};
\node at (4.9,2.2) {$\Longrightarrow$};

\draw[very thick] (6,0) to [out=90,in=195] (7.7,4.5); 
\draw[very thick] (7.7,4.5) to [out=10,in=90] (9.4,0); 
\draw (6.4,3.3) --(9.1,3.5); 
\node at (7.9,4.8) {$p_0$};
\node at (7.9,4.5) {$\bullet$};
\node at (7.7,3.9) {$E'_0$};
\node at (7.7,1.4) {$E'_{per}$};

\node at (10.9,2.7) {$\psi$};
\node at (10.9,2.2) {$\Longrightarrow$};

\draw[very thick] (12,0) to [out=90,in=195] (13.7,4.5); 
\draw[very thick] (13.7,4.5) to [out=10,in=90] (15.4,0); 
\node at (13.9,4.8) {$p_0$};
\node at (13.9,4.5) {$\bullet$};
\node at (13.7,1.4) {$E''_{per}$};
\end{tikzpicture}
\caption{Iterating twice $\psi$ to make $p_{0}$ periodic}
\label{fig:pzeroPeriodic}
\end{figure}
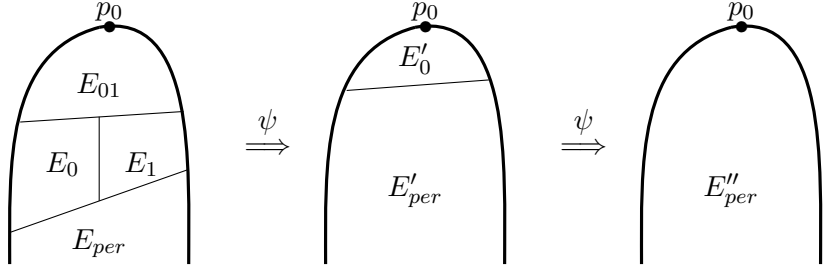
Figure~\ref{fig:pzeroPeriodic} illustrates the main step of the
  proof, the double iteration of $\psi$ to turn $p_{0}$ into a
  periodic point.
Notice that some crucial arguments used in the above proof (starting
from the induction on $\vert e\vert$ itself) make essential use of the
fact that evaluations are order-preserving, so such arguments are not
suitable for modal logics.

\section{A Non-Ultimately Periodic Endomorphism }\label{sec:endomorphisms}

Ruitenburg's Theorem, interpreted over finitely generated free Heyting algebras, says that any endomorphism of such algebras
$$
\mu: \cF_H(x_1, \dots, x_n) \longrightarrow \cF_H(x_1, \dots, x_n)
$$
is ultimately periodic with period 2, in case it fixes all free generators but one. One may ask whether this is a peculiar property of the endomorphisms fixing all free generators but one or whether this can be extended to all endomorphism: we show by a counterexample that there exists endomorphisms of the free algebra over \emph{two} generators which are not periodic.

\begin{figure}
  \centering
  \begin{tikzpicture}
    \OurLadder{3};
    \node at ([xshift=-7mm]origin) {$a$} ;
    \node at ([xshift=-7mm]q0) {$b$} ;
    \node at ([xshift=11mm]q1) {$b \to a$} ;
    \node at ([xshift=-15mm]q2) {$(b \to a) \to b$} ;
    \node at ([xshift=19mm]q3) {$((b \to a) \to b) \to b$} ;
    \node at ([xshift=18mm]q7) {$\da 6 \to (\da 3 \,\vee \da 4)$} ;
\end{tikzpicture}
\caption{The modified \RN ladder}\label{fig:RNmod}
\end{figure}
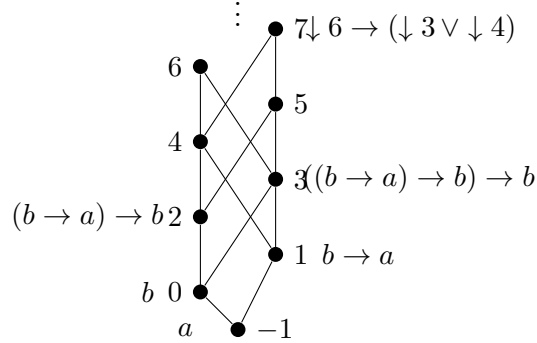

To describe our counterexample we first introduce a variant
$(R, \leq)$ of the Rieger-Nishimura ladder. This is the poset so
described:
\begin{itemize}
 \item[-] $R= \{n\in \Z\mid n\geq -1\}$;
 \item[-] $n\leq m$ iff either $n=-1$ or ($n\geq 0$ and either $n\leq m-2$ or $n=m$).
\end{itemize}
It is not difficult to see that $\leq$ is a reflexive, transitive,
antisymmetric relation (actually, $(R, \leq)$ differs from
Rieger-Nishimura only for the presence of the bottom element $-1$, see Figure~\ref{fig:RNmod}).

Let us consider the Heyting algebra $\cD(R)$
  of
downsets of $(R, \leq)$. We show that this is generated by the two
downsets $a =\{-1\}$ and $b=\{0, -1\}$. To see this, we show
that, for all $n\in R$, the downset $\downarrow n$ belongs to the subalgebra generated by $a,b$.  In fact we have
that:
\begin{itemize}
 \item[-] $\downarrow -1 =a$;
 \item[-] $\downarrow 0=b$;
 \item[-] $\downarrow 1 =b\to a$;
 \item[-] $\downarrow 2 =(b\to a)\to b$;
 \item[-] $\downarrow 3 =((b\to a)\to b)\to b$;
 \item[-] $\downarrow n+4 =\,\downarrow (n+3)\! \to\, \downarrow n \,\vee \downarrow (n+1)$ (for $n\geq 0$).
\end{itemize}

\begin{remark}
  Let $d \in \cD(R)$ be such that $d \neq \emptyset$ and $d \neq R$.
  Then either $d = \downset n$, for some $n \geq -1$, or $d$ has two
  maximal elements $n$ and $n +1$, so $d = \downset n\, \cup \downset(n+1)$. In the
  latter case, both $n + 3$ and $n+4$ are upper bounds of $d$, but
  $n+3$ is such that $\card(\downset(n+3)) = \card(d) + 1$, while
  $\card(\downset(n+4)) = \card(d) + 2$.
  We let therefore:
  \begin{align*}
    \mysup(\downset n) & \eqdef n\,,
    &
    \mysup(\downset n \cup \downset(n + 1)) & \eqdef n+3\,.
  \end{align*}
\end{remark}
\begin{lemma}\label{lem:R}
 $\cD(R)$ is isomorphic to the free Heyting algebra over two generators $a,b$ divided by the congruence generated by
\begin{equation}\label{eq:fpres}
\top = \neg \neg a \wedge (a \to b)~~.
\end{equation}
\end{lemma}

\begin{proof} Let $F$ for the time being be the above mentioned
  finitely presented algebra. Since, within $\cD(R)$, the downset $\{-1\}$ is a least non trivial element (so $\neg \{-1\} = \emptyset$) and
    $\{-1\} \subseteq \{0, -1\}$, it is clear that $\cD(R)$ satisfies the equality~\eqref{eq:fpres} for $a:= \{-1\}$ and $b:= \{0,-1\}$.
    Also, $\cD(R)$
  is generated by the two elements $\{0, -1\}$ and $\{-1\}$; as a consequence, the
  function $q: F\longrightarrow\, \cD(R)$ mapping the
  equivalence class of the free generator $b$ to $\{0, -1\}$ and the
  equivalence class of the free generator $a$ to $\{-1\}$ is a Heyting
  algebras quotient. To show that this quotient map is also injective,
  we use the same technique we adopted for showing that $\cS(h_L)$ is the free algebra over the distributive lattice $\cD(L)$:
  since $F$ embeds into a product of finite Heyting algebras (like any
  finitely presented Heyting algebra, see
  Lemma~\ref{lemma:productOfFinite}),
    it is sufficient to show
  that any morphism $h: F\longrightarrow \,\cD(P)$ (for a
  finite poset $(P,\leq)$) factors through $q$. This follows from the
  following statement on finite Kripke models:
\begin{description}
\item[{\rm (*)} ] for every finite Kripke model over $P$ validating
  $\neg \neg a \wedge (a \to b)$ there is an open map
  $f: (P, \leq)\longrightarrow (R, \leq)$ preserving the evaluation of
  $a,b$.
\end{description}
Property (*) is easily checked by
defining $f(p)$ ($p\in P$) by induction on the height of $p$. In
detail:
\begin{enumerate}
\item[(i)] $f(p)\eqdef-1$ if $p$ forces both $a,b$; 
\item[(ii)] $f(p)\eqdef 0$
  if $p$ forces only $b$;
\end{enumerate}
If $p$ forces  neither $a$ nor $b$, then
\begin{enumerate}
\item[(iii)] if all $p'<p$ force both $a$ and $b$, then
  $f(p)\eqdef 1$; 
\item[(iv)] if all
  $p'<p$ force $b$ and there is $p'<p$ forcing only $b$, then
  $f(p)\eqdef 2$; 
\item[(v)] in all remaining cases, $f(p) \eqdef \mysup \set{ f(p') \mid p'
    < p}$. 
\end{enumerate}
Notice that the above analysis is exhaustive: since $a\to b$ is true
everywhere, there cannot be points forcing only $a$ and not
$b$. Similarly, since $\neg \neg a$ is true everywhere, for every $p$
there must be $p'\leq p$ forcing $a$: this fact is used when checking
that $f$ as defined above is open.
\end{proof}

\begin{lemma}\label{lem:fpnonperiodic}
 There is an endomorphism of $\cD(R)$ which is not ultimately periodic.
\end{lemma}

\begin{proof}
 The endomorphism is the inverse image $f^{-1}$ along the open map
$f:(R, \leq) \longrightarrow (R, \leq)$ so defined:
$$
f(n)\eqdef -1 ~({\rm for}~n<2) \qquad f(n)\eqdef n-2 ~({\rm for}~n\geq 2)~~~.
$$
It is evident that $f^{-1}$ is not ultimately periodic: this comes from the fact that we have $f(n)=n-2$ for all $n\geq 2$. 
\end{proof}

To lift the endomorphism of Lemma~\ref{lem:fpnonperiodic} to the level of the free algebra on two generators, 
we first show that
 $\cD(R)$ is a projective algebra:
 
 \begin{lemma}
   \label{lemma:projective}
 $\cD(R)$ is a projective Heyting algebra.
\end{lemma}
 
 \begin{proof}
   We might use for the proof of this Lemma the general results
   from~\cite{G99}, however we prefer to supply a direct proof. Let
   $F_H(a,b)$ be the free Heyting algebra on two generators and let
   $q: F_H(a,b)\longrightarrow\, \cD(R)$ be the
   homomorphism mapping the free generator $b$ to $\{0, -1\}$ and the
   free generator $a$ to $\{-1\}$: what we have to produce is a
   section of $q$, namely a morphism $s$ in the opposite direction
   such that $q\circ s=id$. Taking into consideration
   Lemma~\ref{lem:R} and turning the existence of $s$ into logical
   terms, what we need is a substitution $\sigma$ such that the
   formulae
 \begin{eqnarray*}
 & \sigma(\neg \neg a \wedge (a \to b))  \\
 & \neg \neg a \wedge (a \to b) \to (a \leftrightarrow \sigma(a)) \\
 & \neg \neg a \wedge (a \to b) \to (b \leftrightarrow \sigma(b))
 \end{eqnarray*}
are provable in intuitionistic logic (here $\sigma(a), \sigma(b)$ must be formulae over the propositional variables $a,b$).
The required substitution in fact exists and can be taken to be  
 $$
a\longmapsto \neg \neg a \to a, \qquad b\longmapsto ((\neg \neg a \to a)\to b)\to b~~
$$
as it can be easily checked.
 \end{proof}
 \begin{theorem}
  There is an endomorphism of the free Heyting algebra on two generators which is not ultimately periodic.
 \end{theorem}

\begin{proof}
  Let $f : \,\cD(R) \ra \cD(R)$ be the morphism defined in
  Lemma~\ref{lem:fpnonperiodic}, so $f^{i+p} \neq f^{i}$ for no
  $i \geq 0$ and $p \geq 1$. Let $q : F_H(a,b) \ra \cD(R)$ and
  $s: \cD(R)\ra F_H(a,b)$ as in the proof of
  Lemma~\ref{lemma:projective}, so $q \circ s = id$. Let
  $g \eqdef s \circ f \circ q$ and suppose that $g^{i + p} =
  g^{i}$. Then $s \circ f^{i+p}\circ q = s \circ f^{i}\circ q$ and, by
  precomposing with $s$ and postcomposing with $q$, $f^{i+p} = f^{i}$, contradiction.
\end{proof}

\section{Bounds for Periods }\label{sec:bounds}

We fix, in this Section, a finite poset $L$ and a natural
transformation $\psi: h_{L} \lora h_{L}$.  We shall pay a particular
attention to the case where $\psi$ is the dual of an endomorphism of a
finitely generated free algebras. This happens exactly when $\psi$ has
a \bindex and when $L$ is of the form
$\langle {\cal P}(\ux), \supseteq\rangle$ for a finite set $\ux$ (so
$\cD(L)$ is a free distributive lattice).

As we saw, $\psi$ might not be ultimately periodic but, on the other
hand, all components $\psi_P$ of $\psi$ are such (because the $P$ are
finite posets). We show that the period of $\psi_P$ can be uniformly
bounded depending on the sole cardinality of $L$
(and not on the cardinality of $P$). 
More precisely, we have
the following statement:
\begin{proposition}\label{prop:uniform}
  Let $\ell$ be the cardinality of $L$.
  For each finite set $P$ and each $v \in h_{L}(P)$, the period of the
  sequence $\set{\psi^{k}_{P}(v) \mid k \geq 0}$ has $\ell!$ as an
  upper bound.
\end{proposition}
The Proposition is an immediate consequence of Lemma~\ref{lem:uniform}
below, for which we need to define a few concepts.

 For a point $p\in P$, we let the
\emph{view set} of $p$ (\wrt $v,\psi$)  be the set
$\set{\psi^{k}(v)(p) \mid k\geq 0}$ and, for $S\subseteq P$, we let
the view set of $S$ (\wrt $v,\psi$)  be the union of
the view sets of the $p\in S$.

Our claim follows from the following:
\begin{lemma}\label{lem:uniform}
  For $v \in h_{L}(P)$, the period of the sequence
  $\set{\psi^{k}(v) \mid k \geq 0}$ has as an upper bound $K!$, where
  $K$ is the cardinality of the view set of $P$.
\end{lemma}
\begin{proof}
  We argue by induction on the height of $P$. If such an height is 1,
  then $P$ contains only the root and the period is bounded by
  $K\leq K!$.
 
  Suppose that the height of $P$ is greater
    than $1$ and let $p$ be the root of $P$. Then, let
  $\ddownset p \eqdef \set{q \in P \mid q < p}$ and let $M$ be the
  cardinality of the view set of $\ddownset p$.
 By the induction hypothesis, for any $q\in \ddownset p$, $M!$ is an
 upper bound for the period of the sequence
 $\set{\psi^{k}(v_{q}) \mid k \geq 0}$. Since the lcm of many copies
 of $M!$ is $M!$, the restriction of $\set{\psi^{k}(v) \mid k \geq 0}$
 to $\ddownset p$ 
 has period $M!$.

 Thus, for $s$ large enough, we have $\psi^{s+M!}(v_q)= \psi^s(v_q)$
 for all $q\in \ddownset p$. Let $a$ be maximal
 (\wrt the partial order of $L$) 
 in the view set of $\ddownset p$ (\wrt $v,\psi^s$). Without loss of
 generality (that is, up to increasing $s$ a bit), we can suppose that
 there is $q_{0}\in \ddownset p$ such that $\psi^s(v)(q_{0})=a$.

 Consider now the set $\set{ \psi^{s+k\cdot M!}(v)(p) \mid k\geq 0}$
 and let $N$ be its cardinality. If, for some $k$,
 $\psi^{s+k\cdot M!}(v)(p)$ belongs to the view set of $\ddownset p$,
 then, for all $q \in \ddownset p$,
 $\psi^{s+k\cdot M!}(v)(q) \leq \psi^{s+k\cdot M!}(v)(p)$, and in
 particular
 $a = \psi^s(v)(q_{0}) = \psi^{s+k\cdot M!}(v)(q_{0}) \leq
 \psi^{s+k\cdot M!}(v)(p)$. It follows that
 $a= \psi^{s+k\cdot M!}(v)(p)$, by the maximality of $a$. Therefore
 the set $\set{ \psi^{s+k\cdot M!}(v)(p) \mid k\geq 0}$ intersects the
 view set of $\ddownset p$ at most in the singleton $\set{a}$ and,
 consequently, we have $M+N-1\leq K$, where $K$ is the cardinality of
 the view set of the whole $P$.  
 We clearly have that $\psi^s$ becomes periodic in at most
 $N\cdot ( M!)$ steps (with period bounded by this number) and
 the claim follows from the inequality $N\cdot (M!)\leq K!$ below.
\end{proof}

\begin{lemma}
 For $M,N\geq 1$, we have $N\cdot ( M!)\leq (M+N-1)!$.
\end{lemma}
\begin{proof}
  The case $N =1$ is obvious, so we suppose that $N > 1$. Since
  $M \geq 1$, $N \leq M + N -1$ 
  and therefore
  \begin{align*}
    N \cdot M! & \leq (M + N-1) M! \\
    & \leq M!(M+1)(M+2) \ldots (M+N-1) = (M+N-1)!\,,
  \end{align*}
  where for the last equality we have used that  $M + N -1 > M$.
\end{proof}

The following result is an immediate consequence of
Proposition~\ref{prop:uniform}:
\begin{proposition}
  \label{prop:periodLocallyFinite}
  Let a finitely generated free \Ha
  homomorphism
  $\mu: \cF_H(x_1, \dots, x_n) \longrightarrow \cF_H(x_1, \dots, x_n)$
  be ultimately periodic; then its period is bounded by $2^n!$.
\end{proposition}

\begin{remark}  
  Let us point out that the bound given in
  Proposition~\ref{prop:uniform} strictly depends on
  $h_{L}(P)$ being a set of monotone functions.  Observe that, when
  the bound has been constructed in the proof of
  Lemma~\ref{lem:uniform}, the function
  $\psi_{P} : h_{L}(P) \ra h_{L}(P)$ has been decomposed as
  $\psi_{P}(\vec{y},x) = (g(\vec{y}),f(\vec{y},x))$, where $\vec{y}$
  is a vector of elements of $L$ indexed by $\ddownset p$ and
  $x \in L$; moreover $\psi$ is applied to pairs $(\vec{y},x)$
  such that $y \leq x$ for each $y \in \vec{y}$.
  If we give away the latter constraint on the order, it is easy to
  see that the bound does not hold anymore. This happens, even when
  $\psi$ is recursively defined on the height of $P$ (so that all of
  its restrictions $\psi_{p}$ are of the form
  $\langle g \circ\pi_{1},f\rangle$ for some
  $g : \ddownset p \ra \ddownset p$ and for some
  $f : L^{\ddownset p} \times L \rto L$).
  Consider the following example.  Let $P$ be the chain
  $\set{1,\ldots ,n}$ and let $L = \set{0,1}$, so we can identify \emph{arbitrary}
  functions from $P$ to $L$ with words on the alphabet $\set{0,1}$ of
  length $n$.  For $x \in \set{0,1}$, let $\psi_{1}(x) \eqdef 1- x$.
  Suppose that, for $i < n$, we have defined
  $\psi_{i} : \set{0,1}^{i} \rto \set{0,1}^{i}$ so that $\psi_{i}$ is a
  bijection of period/order $2^{i}$. We can list then
  $\set{0,1}^{i} = \set{w_{0},\ldots ,w_{2^{i}-1}}$ with
  $w_{j} = \psi_{i}^{j}(0,\ldots ,0)$, $j = 0,\ldots , 2^{i}-1$. Define then
  \begin{align*}
    \psi_{i+1}(w_{j},x) & \eqdef
    \begin{cases}
      (h_{i}(w_{j}),x)\,, & j < 2^{i} -1\,, \\
      (h_{i}(w_{j}),1- x)\,, & j = 2^{i}-1\,.
    \end{cases}
  \end{align*}
  This recursive construction yields $\psi_{n}$ of period $2^{n}$ and,
  in particular, the factorial bound $2 !$ for the period does not
  apply.  
\end{remark}

A subvariety $\bf V$ of Heyting algebras is said to be \emph{locally
  finite} iff the finitely generated free $\bf V$-algebras are all
finite (we shall indicate with $\cF_{\bf V}(x_1, \dots, x_n)$ the free
$\bf V$-algebra on the generators $x_1, \dots, x_n$. Obviously, a
locally finite subvariety is also finitely approximable, hence
Theorem~\ref{thm:Vduality} applies to it. Since all results in this
section trivially apply also to finitely approximable varieties and
since the endomorphisms between finitely generated free
$\bf V$-algebras are ultimately periodic (by the finiteness of these
algebras), we obtain:

\begin{theorem} Let $\bf V$ be  a locally finite variety of Heyting algebras. Every
  endomorphism
  $\mu: \cF_{\bf V}(x_1, \dots, x_n) \longrightarrow \cF_{\bf V}(x_1, \dots, x_n)$
  is ultimately periodic and its period is bounded by $2^n!$.
\end{theorem}

The interesting point is that the above bound is uniform with respect
to all varieties $\bf V$.  However it is not in general tight, as we
shall remark in the final section by considering as ${\bf V}$ the
variety of Boolean algebras.

\section{Conclusions and Open Problems}

\RT exhibits a particular finitistic behaviour of one-variable
substitutions in the \IPC.  Willing to provide a semantical proof of
this theorem, we have studied more general substitutions, which,
algebraically, can be identified with endomorphisms of finitely
generated free \Ha{s}.  The proof of \RT as well as some additional
remarks on periods of iterated substitutions have been achieved using
the semantical apparatus given by the sheaf theoretic duality for
finitely presented \Ha{s} \cite{GhilardiZawadowski2011}.

Using these semantical tools, sheaf duality and bounded bisimulations,
we found upper bounds for the index and the period of sequences of
iterated substitutions.  The bounds so found are not optimal.  For
example, the proof of Theorem~\ref{thm:main} yields a bound for
the index which is non elementary as a function of
the implication degree of an IPC formula.  On the other hand, the
bound that can be
  extracted from the syntactic computations in~\cite{Ruitenburg84} is
linear \wrt the implicational degree and the number of propositional
variables of a formula. The syntactic computations in~\cite{fossacs}
for fixpoints convergence also yield tighter bounds.

While the semantical approach has been successful for providing a
proof of \RT, it remains open whether similar approaches can yield
finer bounds.

Other open problems arise from inspecting the results presented in this paper.
Firstly, although we were able to show that periodicity fails for two-variable substitutions, it is still an open problem to characterize 
or to decide periodicity for abitrary substitutions in IPC (the only sufficient condition known is the one supplied by \RT - namely the fact that all-but-one variables are fixed).

Secondly, concerning the upper bounds we found, notice that Proposition~\ref{prop:periodLocallyFinite}
provides bounds for \emph{periods} of free \Ha endomorphisms in
locally finite subvarieties; being able to bound their \emph{indexes}
might also be interesting. In particular, it is not clear whether
indexes are sensitive to the number of generators of a free algebra
or, similarly to what happen for fixpoint approximants in some lattice
varieties, see \cite{ramics2014}, they are uniform in a fixed
variety. Corollary~\ref{cor:indexFromHeight} can be used to argue that
this is the case in varieties of \Ha{s} of bounded height---see
e.g. the varieties 
$\mathbf{bd_n}$ in \cite[Prop. 2.38]{CZ}---yet there are locally
finite varieties of \Ha{s}---notably, the variety of G\"odel/Dummet
algebras---whose Kripke models might be of unbounded height.

Coming back to the period, let us also notice that the upper bound
provided by Proposition~\ref{prop:periodLocallyFinite} is not in
general tight. To see why, consider free Boolean algebras: a morphism
$f : \cF_{B}(x_{1},\ldots ,x_{n}) \ra \cF_{B}(x_{1},\ldots ,x_{n})$
corresponds, via duality, to a function $f : \two^{n} \rto \two^{n}$.
Now, estimating an upper bound for the periods of functions from the
set $[k] \eqdef \set{1,\ldots ,k}$ (where in our case $k = 2^{n}$) to
itself can be reduced to estimating an upper bound for the period (or
order) of permutations of $[k]$. Indeed, if $i$ and $p$ are such that
$f^{i+p} = f^{i}$, then the restriction of $f$ to $f^{i}([k])$ is a
permutation of $f^{i}([k])$ which can be extended to a full
permutation of the set $[k]$ of equal period $p$.  Now, an upper bound
for all these periods is $\lcm(1,\ldots ,k)$ for which we have
$2^{k-1} \leq \lcm(1,\ldots ,k) \leq 3^{k}$ \cite{Farhi} and,
asymptotically, $\lcm(1,\ldots ,k) \sim e^{k}$ (by the prime number
theorem).
On the other hand, using Stirling approximation,
$k! \sim \sqrt{2\pi k}(\frac{k}{e})^{k}$.
It is an open problem whether the bound $2^{n}!$ can be made tighter
by considering locally finite varieties of \Ha{s} other than Boolean
algebras; it is not clear either how the bound can vary below $2^{n}!$
depending on the locally finite subvariety $\bf V$.

Finally, most of the techniques used here for \Ha{s} are also the
tools for studying modal logics in
\cite{GhilardiZawadowski2011}. While we can expect that periodicity
phenomena of substitutions do not arise for the basic modal logic
$\mathbf{K}$, they surely do for locally tabular modal logics.
Considering also the numerous results on definability of fixpoints,
see e.g. \cite{Sambin1976,AlberucciFacchini2009b}, these phenomena are
likely to appear in other subsystems of modal logics.  As far as we
know, investigation of periodicity phenomena in modal logics is a
research direction which has not yet been explored and where the
bounded bisimulation methods might prove their strength once more.

\ifams
\bibliographystyle{abbrv}
\fi
\ifmscs
\bibliographystyle{apalike}
\fi
\bibliography{biblio,biblio2}

\def\cprime{$'$}
\begin{thebibliography}{10}

\bibitem{AlberucciFacchini2009b}
L.~Alberucci and A.~Facchini.
\newblock The modal {$\mu$}-calculus hierarchy over restricted classes of
  transition systems.
\newblock {\em J. Symbolic Logic}, 74(4):1367--1400, 2009.

\bibitem{birkhoff1937}
G.~Birkhoff.
\newblock Rings of sets.
\newblock {\em Duke Math. J.}, 3(3):443--454, 09 1937.

\bibitem{CZ}
A.~Chagrov and M.~Zakharyaschev.
\newblock {\em Modal logic}, volume~35 of {\em Oxford Logic Guides}.
\newblock The Clarendon Press, Oxford University Press, New York, 1997.
\newblock Oxford Science Publications.

\bibitem{DP}
B.~A. Davey and H.~A. Priestley.
\newblock {\em Introduction to lattices and order}.
\newblock Cambridge University Press, New York, second edition, 2002.

\bibitem{Esa74}
L.~Esakia.
\newblock Topological {K}ripke models.
\newblock {\em Soviet Math. Dokl.}, 15:147--151, 1974.

\bibitem{Farhi}
B.~Farhi.
\newblock An identity involving the least common multiple of binomial
  coefficients and its application.
\newblock {\em Amer. Math. Monthly}, 116(9):836--839, 2009.

\bibitem{fine}
K.~Fine.
\newblock Logics containing {${\rm K}4$}. {II}.
\newblock {\em J. Symbolic Logic}, 50(3):619--651, 1985.

\bibitem{ramics2014}
S.~Frittella and L.~Santocanale.
\newblock Fixed-point theory in the varieties $\mathcal{D}_{n}$.
\newblock In P.~H{\"o}fner, P.~Jipsen, W.~Kahl, and M.~E. M{\"u}ller, editors,
  {\em RAMICS}, volume 8428 of {\em Lecture Notes in Computer Science}, pages
  446--462. Springer, 2014.

\bibitem{unpro}
S.~Ghilardi.
\newblock Unification through projectivity.
\newblock {\em J. Logic Comput.}, 7(6):733--752, 1997.

\bibitem{G99}
S.~Ghilardi.
\newblock Unification in intuitionistic logic.
\newblock {\em J. Symbolic Logic}, 64(2):859--880, 1999.

\bibitem{um}
S.~Ghilardi.
\newblock Best solving modal equations.
\newblock {\em Ann. Pure Appl. Logic}, 102(3):183--198, 2000.

\bibitem{unapal}
S.~Ghilardi.
\newblock Unification, finite duality and projectivity in varieties of
  {H}eyting algebras.
\newblock {\em Ann. Pure Appl. Logic}, 127(1-3):99--115, 2004.
\newblock Provinces of logic determined.

\bibitem{fossacs}
S.~Ghilardi, M.~J. Gouveia, and L.~Santocanale.
\newblock Fixed-point elimination in the intuitionistic propositional calculus.
\newblock In {\em Foundations of Software Science and Computation Structures,
  {FOSSACS} 2016, Proceedings}, pages 126--141, 2016.

\bibitem{aiml18}
S.~Ghilardi and L.~Santocanale.
\newblock Ruitenburg's theorem via duality and bounded bisimulations.
\newblock In {\em Advances in Modal Logic, {AiML} 2018, Proceedings}, pages
  277--290, 2018.

\bibitem{GhilardiZawadowski2011}
S.~Ghilardi and M.~Zawadowski.
\newblock {\em Sheaves, Games, and Model Completions: A Categorical Approach to
  Nonclassical Propositional Logics}.
\newblock Springer Publishing Company, Incorporated, 1st edition, 2011.

\bibitem{GhilardiZawadowski97}
S.~Ghilardi and M.~W. Zawadowski.
\newblock Model completions, r-{H}eyting categories.
\newblock {\em Ann. Pure Appl. Logic}, 88(1):27--46, 1997.

\bibitem{goguen}
J.~A. Goguen.
\newblock What is unification? {A} categorical view of substitution, equation
  and solution.
\newblock In {\em Resolution of equations in algebraic structures, {V}ol. 1},
  pages 217--261. Academic Press, Boston, MA, 1989.

\bibitem{Mardaev07}
S.~Mardaev.
\newblock Definable fixed points in modal and temporal logics : {A} survey.
\newblock {\em Journal of Applied Non-Classical Logics}, 17(3):317--346, 2007.

\bibitem{Mardaev1993}
S.~I. Mardaev.
\newblock Least fixed points in \grzl and in the intuitionistic propositional
  logic.
\newblock {\em Algebra and Logic}, 32(5):279--288, 1993.

\bibitem{Ruitenburg84}
W.~Ruitenburg.
\newblock On the period of sequences $(a^n (p))$ in intuitionistic
  propositional calculus.
\newblock {\em The Journal of Symbolic Logic}, 49(3):892--899, Sept. 1984.

\bibitem{Sambin1976}
G.~Sambin.
\newblock An effective fixed-point theorem in intuitionistic diagonalizable
  algebras.
\newblock {\em Studia Logica}, 35(4):345--361, 1976.
\newblock The algebraization of the theories which express Theor, IX.

\bibitem{shavrukov}
V.~Y. Shavrukov.
\newblock Subalgebras of diagonalizable algebras of theories containing
  arithmetic.
\newblock {\em Dissertationes Math. (Rozprawy Mat.)}, 323:82, 1993.

\bibitem{visser}
A.~Visser.
\newblock Uniform interpolation and layered bisimulation.
\newblock In {\em G\"odel '96 ({B}rno, 1996)}, volume~6 of {\em Lecture Notes
  Logic}, pages 139--164. Springer, Berlin, 1996.

\end{thebibliography}

\end{document}